\documentstyle{amsppt}
\input amstex
\hsize= 5.5in
\NoBlackBoxes
\input epsf
\topmatter 
\title A non-quasiconvexity embedding theorem for hyperbolic groups.\endtitle
\rightheadtext {Non-Quasiconvexity Embedding Theorem}
\author Ilya Kapovich
\endauthor
\address 
Department of Mathematics, Hill Center, Busch Campus,
Rutgers University, Piscataway, NJ 08854
\endaddress

\email ilia\@math.rutgers.edu\endemail
\abstract 
We show that if $G$ is a non-elementary torsion-free word hyperbolic
group then there exists another word hyperbolic group $G^*$, such that
$G$ is a subgroup of $G^*$ but $G$ is not quasiconvex in $G^*$.
\endabstract

\endtopmatter

\document

\head 1.Introduction. \endhead

Word hyperbolic groups introduced by M.Gromov in \cite{Gr} have played a very
important role in the recent progress of geometric and combinatorial
group theory. It has turned out that many groups arising in a
traditionally geometric context, such as fundamental groups of closed
compact manifolds admitting a metric of strictly negative curvature,
have the property of being word hyperbolic. On the other hand, most
finitely presented groups (in a certain probabilistic sense), also
belong to the class of word hyperbolic groups. Moreover, word
hyperbolic groups possess a number of very good algorithmic and
combinatorial properties which do not hold for the class of finitely
presented groups in general. For example, the word problem and the
conjugacy problem are solvable in any word hyperbolic group \cite{Gr, 7.4.B}. Even
more surprisingly, the isomorphism problem is solvable in the class of
torsion-free word hyperbolic groups \cite{Sel}.
We will provide some elementary facts about hyperbolic groups in
Section 2. For a more detailed discussion the reader is referred to
\cite{Gr}, \cite{ABC}, \cite{CDP}, \cite{GH}.

\medskip
A particularly important and interesting class of subgroups of word
hyperbolic groups are the so-called {\it quasiconvex subgroups}.
Roughly speaking, a finitely generated subgroup $H$ of a word
hyperbolic group $G$ is quasiconvex in $G$ if for any word metric
$d_G$ on $G$ and any word metric $d_H$ on $H$ the inclusion map
$i:(H,d_H)\longrightarrow (G,d_G)$ is a bi-Lipschitz embedding.
(see Section 2 for a more careful definition).
\medskip
Quasiconvex subgroups of word hyperbolic groups have many good
properties. For example, a quasiconvex subgroup of a word hyperbolic
group is always finitely presented and itself word hyperbolic. The
intersection of any two quasiconvex subgroups is always quasiconvex
(and so is finitely generated and finitely presented). Quasiconvex
subgroups are of particularly great importance when one studies
amalgamated free products and HNN-extensions of hyperbolic groups
\cite{Gr}, \cite{BF1}, \cite{KM}, \cite{BGSS}, \cite{Gi}, \cite{Pa}.
\medskip
From a geometric standpoint, quasiconvex subgroups of word hyperbolic
groups are closely related to geometrically finite subgroups of
classical hyperbolic groups. More precisely, suppose $G$ is a
geometrically finite group of isometries of ${\Bbb H}^n$ without
parabolics. Then $G$ is word hyperbolic and a subgroup $H$ of $G$ is
quasiconvex in $G$ if and only if $H$ is geometrically finite \cite{Swa}.
Some of the properties of quasiconvex subgroups will be given in
Section 3. For a more detailed discussion on quasiconvexity the reader
is referred to \cite{Gr}, \cite{GH}, \cite{KS}, \cite{MT}.
\medskip
One of the most fundamental results pertaining to quasiconvexity is
the following theorem of M.Gromov \cite{Gr, 8.1.D} (see also \cite{ABC,
Proposition 3.2}).

\proclaim{Proposition 1.1} Let $G$ be a word hyperbolic group and let
$C$ be an infinite cyclic subgroup of $G$. Then $C$ is quasiconvex in
$G$.
\endproclaim
Thus for an infinite cyclic group quasiconvexity is an absolute
property and does not depend on the embedding in an ambient hyperbolic
group.
This leads us to the following definition.
\definition{Definition 1.2} Let $G$ be a word hyperbolic group. We say
that $G$ is {\it absolutely quasiconvex} if for any word hyperbolic
group $G^*$, containing $G$ as a subgroup, $G$ is quasiconvex in
$G^*$.
\enddefinition
It is natural to ask, then, whether there are any absolutely
quasiconvex groups other than the infinite cyclic group. Perhaps
surprisingly, it turns out that the answer to this question is no, at
least if we restrict ourselves to the class of torsion-free word
hyperbolic groups.
The main result of the present paper is the following
\proclaim{Theorem A} Let $G$ be a non-elementary (that is not
virtually cyclic) torsion-free word hyperbolic group. Then there
exists a word hyperbolic group $G^*$ such that $G$ is a subgroup of
$G^*$ and $G$ is not quasiconvex in $G^*$.
\endproclaim
Since any non-trivial torsion-free virtually cyclic group is in fact
infinite cyclic, this result immediately implies
\proclaim{Theorem B} Let $G$ be a non-trivial torsion-free absolutely
quasiconvex group. Then $G$ is infinite cyclic.
\endproclaim
We would like to note that up to this point there were very few known
examples of finitely generated subgroups of word hyperbolic groups
that are not quasiconvex. 
\smallskip
The first example of this kind is provided in a remarkable work of E.Rips \cite{Ri}.
Given any finitely presented group $Q$, Rips constructs a finitely
generated small cancellation $C(7)$-group $G$ and a short exact
sequence
$$1\rightarrow K\rightarrow G\rightarrow Q\rightarrow 1$$
where $K$ is a two-generated group. The group $G$ is word hyperbolic
by the basic results of small cancellation theory.
It was shown in \cite{ABC} that if a normal subgroup of a word hyperbolic
group is quasiconvex then it is either finite or has finite index. 
Suppose now that in the Rips' construction the group $Q$ is chosen to
be a non-hyperbolic infinite group (e. g. ${\Bbb Z}\times {\Bbb Z}$). Since $G$ is word hyperbolic but $G/K$ is not, we conclude that
$K$ is infinite. Thus $K$ is an infinite two-generated subgroup of
infinite index in $G$ which is normal in $G$. Therefore $K$ is not
quasiconvex.Therefore in Rips's example $H$ is not quasiconvex. It was
later noticed in \cite{Sho} and in \cite{BMS}, that $K$ is not even finitely presentable.
\smallskip
  If $M$ is a 3-manifold obtained from a cylinder over a closed oriented surface $S$ of genus at least two by gluing upper and lower boundaries of this cylinder along a pseudo-anosov homeomorphism of $S$ then the fundamental group $G$ of the resulting manifold $M$ is word hyperbolic and it contains a subgroup $H$ isomorphic to the fundamental group of $S$ which is not quasiconvex in $G$. This follows from the result of Thurston \cite{Th} which asserts that in this situation $M$ admits a metric of constant negative curvature. Therefore (see \cite{Gr}) $G$ is word hyperbolic.
Besides $M$ fibers over a circle with a fiber $S$ and thus there is a short exact sequence $1\rightarrow \pi_1(S)\rightarrow G\rightarrow {\Bbb Z}\rightarrow 1$.  Hence by the above mentioned result of \cite{ABC} $\pi_1(S)$ is not quasiconvex in $G$.

Later M.Bestvina and M.Feign \cite{BF1} showed that there is an analog of
the previous example for free groups, that is there is an HNN-extension of a free group by a
suitable automorphism which is word hyperbolic. Indeed, their result shows that if $F$ is a finitely generated noncyclic free group then there is a word hyperbolic group $G$ and a short exact sequence $1\rightarrow F\rightarrow G\rightarrow {\Bbb Z}\rightarrow 1$ and, of course, $F$ is not quasiconvex in $G$.

The existence of a hyperbolic 3-manifold fibering over a circle (see
the example of W.Thurston above) provided a basis for constructing
other examples of nonquasiconvex subgroups of hyperbolic groups. In
\cite{BM} G.Mess and B.Bowditch showed that there is discrete cocompact
group $G$ of isometries of ${\Bbb H}^4$ containing a subgroup $H$
which is finitely generated but not finitely presentable and therefore
not quasiconvex. L.Potyagailo constructed in \cite{Po} a geometrically
finite subgroup $G$ of $SO(4,1)$ without parabolics which contains a
finitely generated subgroup $H$ which is not finitely presentable,
contains infinitely many conjugacy classes of finite subgroups and
such that $G/H={\Bbb Z}$;  (clearly this $H$ is not quasiconvex).
\smallskip
Recently N.Brady \cite{Br} constructed an example of a word
hyperbolic group $G$ which possesses a finitely presented subgroup $H$
such that $H$ itself is not word hyperbolic. As we mentioned before,
quasiconvex subgroups of word hyperbolic groups are word
hyperbolic. Therefore the subgroup $H$ is obviously not quasiconvex in $G$.
\medskip
Thus Theorem A provides a large class of new examples of finitely presented
non-quasiconvex subgroups of hyperbolic groups. Theorem B shows that
there are no rigid torsion-free groups (except for cyclic groups) in
the sense of always being quasiconvex. 
\medskip
We already mentioned the parallel between quasiconvexity and
geometric finiteness. It is interesting to note that there are also
very few examples of non-geometrically finite groups acting discretely
on ${\Bbb H}^n$. There is some reason to believe that, in contrast to
Theorem B, there may exist some ``absolutely geometrically finite''
groups that are not virtually cyclic. Consider, for example, a group
$G$ which is the fundamental group of a closed hyperbolic 3-manifold
fibering over a circle. It would be interesting to investigate whether
$G$ admits a discrete non-geometrically finite action in some
higher-dimensional hyperbolic space ${\Bbb H}^n$, $n>3$.

\head 2.Preliminary facts and definitions. \endhead
\definition {Definition 2.1} If $(X,d)$ is a metric space, then a map
$f\colon [a,b]\rightarrow X$ is called a {\it geodesic segment} if for
any $[a_1,b_1]\subset [a,b]$ $$|a_1-b_1|=d(f(a_1),f(b_1)).$$ 
(Sometimes, by abuse of notation, we will identify such a map $f$ with
its image and denote it by $[f(a), f(b)]$.)
\smallskip
A metric
space $(X,d)$ is said to be {\it geodesic} if any two points in $X$ can be
joined by a geodesic segment.
\smallskip
A geodesic space $(X,d)$ is called $\delta$-{\it hyperbolic} if for each triangle $\Delta$ with geodesic sides in $X$ and for any point $x$ on one of the sides of $\Delta$ there is a point $y$ on one of the two other sides such that $d(x,y)\le \delta$.
\enddefinition
If $G$ is a group and $X=S\cup S^{-1}$ is a finite generating set for $G$ then a {\it Cayley graph} $\Gamma (G,X)$ of $G$ is an oriented labeled graph with $\lbrace g| g\in G\rbrace$ as a set of vertices and an oriented edge $e=(g,ga)$ labeled by $a$ for each $g\in G$, $a\in X$. It is not hard to see that $\Gamma (G,X)$ is a connected locally finite graph. If we put each edge to be isometric to a unit interval, we can define the length of an edge-path in $\Gamma (G,X)$. 
Now for any vertices $g,h$ of $\Gamma (G,X)$ put $d_X(g,h)$ to be the minimal length of an edge-path connecting $g$ to $h$. Then $d_X$ is a metric on $G$ which can be naturally extended to a metric (which we also denote $d_X$) on $\Gamma (G,X)$. The metric $d_X$ is called a {\it word metric} corresponding to $X$.
It is easy to see that $(\Gamma (G,X),d_X)$ is a geodesic metric space.
We denote the set of all words over $X$ (that is the free monoid on
$X$) by $X^{\ast}$ and the element of $G$ represented by a word $w\in
X^{\ast}$ by $\overline w$. For every $g\in G$ we denote
$l_X(g)=d_X(g,1)$ and call $l_X(g)$ the {\it word length} of $G$ with
respect to $X$. Note that for any $g\in G$ the word length $l_X(g)$ is the minimal
number $n\ge 0$ such that $g$ can be expressed as a product 
$$g=x_1\dots x_n,$$ where all $x_i\in X$.
\definition {Proposition-Definition 2.2} (see [ABC, Section 2] for proof)
Let $G$ be a finitely generated group. Then the following conditions are equivalent:
\roster
\item "(1)" for some finite generating set $X=S\cup S^{-1}$ of $G$ and for some $\delta\ge 0$ the Cayley graph $\Gamma (G,X)$ is $\delta$-hyperbolic;
\item "(2)" for any finite generating set $X=S\cup S^{-1}$ of $G$ there is $\delta\ge 0$ such that the Cayley graph $\Gamma (G,X)$ is $\delta$-hyperbolic.
\endroster
If any of these conditions is satisfied, the group $G$ is called {\it word hyperbolic}. 
\enddefinition
It turns out that the class of word hyperbolic groups is very large.
\proclaim{Proposition 2.3}
\roster
\item Finitely generated free groups are word hyperbolic \cite{Gr,
Section 3}.
\item Finite groups are word hyperbolic \cite{Gr, Section 3}.
\item 
Any finitely generated group satisfying one of the  $C(7)$,
$C(5)-T(4)$, $C(4)-T(5)$, $C(3)-T(7)$ small
cancellation conditions is word
hyperbolic \cite{GS1}.
\item If $A$ and $B$ are groups such that $A\cap B$ has finite index
in both $A$ and $B$ then $A$ is word hyperbolic if and only if $B$ is
word hyperbolic.
\item If we have a short exact sequence

$$1\longrightarrow N\longrightarrow G\longrightarrow Q\longrightarrow
1$$ 
where $N$ is finite, then $G$ is word hyperbolic if and only if $Q$ is
word hyperbolic.
\item A free product of two word hyperbolic groups is word hyperbolic
\cite{BGSS, Theorem H}.
\item A free factor of a word hyperbolic group is word hyperbolic \cite{BGSS, Theorem H}.
\item Let $X$ be a geodesic metric space that is $\delta$-hyperbolic
for some $\delta \ge 0$. Let $G$ be a group acting on $X$ by
isometries so that the quotient space $X/G$ is compact and the action
is properly discontinuous, that is for any compact subset $K$ of $X$
the set
$$\{g\in G | K\cap gK\ne \emptyset\}$$ is finite.
Then $G$ is word hyperbolic \cite {CDP,
Corollary 4.5}.
\item A finitely presented group $G$ is word hyperbolic if and only if
$G$ satisfies linear isoperimetric inequality \cite{ABC, Section 2}.
\item A finitely presented group $G$ is word hyperbolic if and only if
$G$ possesses a finite Dehn presentation \cite{ABC, Theorem 2.12,
Theorem 2.5}.
\item Let $G$ be the fundamental group of a closed Riemannian manifold
all of whose sectional curvatures are negative. Then $G$ is word
hyperbolic \cite{Gr, Section 2}. 
\item Let $G$ be a geometrically finite group of isometries of ${\Bbb
H}^n$ without parabolics. Then $G$ is word hyperbolic \cite{Swa}.
\item Let $G$ be a uniform lattice in a rank one real semi-simple Lie group. Then
$G$ is word hyperbolic \cite{Gr, Section 2}.
\item Let $G$ be a finitely generated Coxeter group which does not contain free abelian
subgroups of rank 2. Then $G$ is word hyperbolic \cite{Mou}.
\endroster
\endproclaim
We will list some of the good properties of word hyperbolic groups in
the following statement.  
\proclaim{Proposition 2.4} Let $G$ be a word hyperbolic group. Then
the following holds.
\roster
\item The group $G$ is finitely generated and finitely presented
\cite{Gr, 2.2.A}, \cite{ABC, Proposition 2.10}
\item The group $G$ has solvable word problem and solvable conjugacy
problem \cite{ABC, Theorem 2.18}, \cite{Gr,7.4.B}.
\item The group $G$ has only finitely many conjugacy classes of
elements of finite order \cite{Gr, 2.2.B}, \cite{GH, Section 4,
Proposition 13}.
\item Any subgroup of $G$ is either virtually cyclic or contains a
free subgroup of rank 2 \cite{GH, Section 8, Theorem 37}.
\item The group $G$ is automatic \cite{ECH, Theorem 3.4.5}.
\item The group $G$ does not contain free abelian subgroups of rank 2
and , more generally, the group $G$ does not contain any
Baumslag-Solitar subgroups
$$B(m,n)=<a,t | t^{-1}a^n t =a^m>$$
where $n\ne 0, m\ne 0$ \cite{GS2, Theorem 5.1, Theorem 6.7}.
\item The group $G$ has the property $FP_{\infty}$ \cite{ABC,
Corollary 4.12(2)}.
\item The rational homology $H_*(G,{\Bbb Q})$ and cohomology
$H^*(G,{\Bbb Q})$ groups of $G$ are finite dimensional \cite{ABC,
Corollary 4.12(5)}.
\item If $G$ is torsion-free then $G$ has finite cohomological
dimension. If $G$ is virtually torsion-free, then $G$ has finite
virtual cohomological dimension \cite{ABC,
Corollary 4.12(3)(4)}. 
\item The isomorphism problem is algorithmically solvable in the class
of torsion-free word hyperbolic groups \cite{Sel}.
\endroster
\endproclaim
\medskip
\definition{Definition 2.5} Let $(X,d)$ be a metric space and let
$\lambda>0, C>0$. Let $[a,b]\subset {\Bbb R}$ be a closed interval.
We say that a naturally parameterized path 
$$\alpha: [a,b]\longrightarrow X$$
is a $(\lambda, C)$-quasigeodesic if for any $s,t\in [a, b]$ we have
$$|t-s|\le \lambda d(\alpha(s), \alpha(t))+C.$$
\enddefinition
Let $G$ be a finitely generating group with a finite generating set
$X=S\cup S^{-1}$. Then any word $w$ over $X$ corresponds to a path $p(w)$ in $\Gamma(G,X)$
from $1$ to $\overline w$ of length $l(w)$. 
We say that $w$ is {\it geodesic} with respect to $d_X$ if the path
$p(w)$ is geodesic in $\Gamma(G,X)$, that is if $l(w)=l_X(\overline
w)$.
Similarly, we say that a word $w$ over $X$ in $(\lambda,
C)$-{\it quasigeodesic} with respect to $d_X$ if the path $p(w)$ is $(\lambda,
C)$-quasigeodesic in $\Gamma(G,X)$.
\medskip
A very important property of quasigeodesics in hyperbolic metric
spaces is the following
\proclaim{Proposition 2.6}\cite{ABC, Proposition 3.3}. Let $(X,d)$ be a geodesic metric space that
is $\delta$-hyperbolic for some $\delta\ge 0$. Let $\lambda >0$, $C>0$
be some positive numbers. Then there exists a constant $\epsilon >0$
with the following properties.
\smallskip
Suppose $\alpha:[a,b]\longrightarrow X$ and
$\beta:[a',b']\longrightarrow X$ are $(\lambda,
C)$-quasigeodesics in $X$ such that $\alpha(a)=\beta(a')$ and
$\alpha(b)=\beta(b')$. 
\smallskip
Then for every $x\in im(\alpha)$ there is $y\in
im(\beta)$ such that $$d(x,y)\le \epsilon$$
and for every $y\in im(\beta)$ there is $x\in im(\alpha)$ such that 
$$d(x,y)\le \epsilon.$$
(That is the sets $im(\alpha)$ and $im(\beta)$ are
$\epsilon$-Hausdorff close.)
\endproclaim

For the rest of this article, unless specified otherwise, we will
assume that all finite generating sets of groups that we consider
are closed under taking inverses.

\head 3.Quasiconvex subgroups of word hyperbolic groups. \endhead

In general, a subset $Y$ of a geodesic metric space $(X,d)$ is termed 
$\epsilon$-{\it quasiconvex} in $X$ if for any $y_1, y_2\in Y$ any geodesic 
$[y_1,y_2]$ 
in $X$, joining points $y_1$ and $y_2$, is contained in the
$\epsilon$-neighborhood of $Y$.

\definition {Proposition-Definition 3.1} (see [ABC] for proof)
Let $G$ be a word hyperbolic group and $A$ be a subgroup of $G$. Then the following conditions are equivalent:
\roster
\item  for {\bf some} finite generating set $X=S\cup S^{-1}$ of $G$ there is
an $\epsilon \ge 0$ such that $A$ is an $\epsilon$-quasiconvex subset
of $\Gamma (G,X)$;
\item  for {\bf any} finite generating set $X=S\cup S^{-1}$ of $G$ there is
an $\epsilon \ge 0$ such that $A$ is an $\epsilon$-quasiconvex subset
of $\Gamma (G,X)$;
\item the group $A$ is finitely generated and for {\bf some} finite
generating set $Y=T\cup T^{-1}$ of $A$ and for {\bf some} finite generating set $X=S\cup S^{-1}$ of $G$ there is a constant $C>0$ such that for any $a\in A$
$$ d_X(a,1)/C\le d_Y(a,1)\le Cd_X(a,1);$$ 
\item the group $A$ is finitely generated and for {\bf any} finite
generating set $Y=T\cup T^{-1}$ of $A$ and for {\bf any} finite generating set $X=S\cup S^{-1}$ of $G$ there is a constant $C>0$ such that for any $a\in A$
$$ d_X(a,1)/C\le d_Y(a,1)\le Cd_X(a,1).$$ 
\endroster
If any of these conditions is satisfied then $A$ is called a {\it quasicovex subgroup} of $G$.
\enddefinition
We summarize the basic properties of quasiconvex subgroups of
word hyperbolic groups in the following statement.
\proclaim{Proposition 3.2}

Let $G$ be a word hyperbolic group.
\roster

\item  \cite{ABC, Lemma 3.8} If $A$ is a quasiconvex subgroup of $G$ 
then $A$ is finitely generated, finitely presented and word hyperbolic.

\item  \cite{Sho} If $A_1, \dots, A_k$ are quasiconvex
subgroups of $G$,
then $A=A_1\cap\dots\cap A_k$ is quasiconvex in $G$
 
\item  If $A$, $B$ are subgroups of $G$  
such that $A$ has finite index in $B$,
then $A$ is quasiconvex in $G$ 
if and only if $B$ is quasiconvex in $G$.

\item  \cite{ABC, Corollary 3.4} If $C$ is a virtually cyclic subgroup of 
$G$
 then $C$ is quasiconvex in $G$.

\item \cite{GS2} If $H$ is a subgroup of finite index in $G$ then $H$ is
quasiconvex in $G$. 

\item If $H$ is quasiconvex in $G$ and $K$ is a subgroup of $H$, $K\le
H\le G$, then $K$ is quasiconvex in $G$ if and only if $K$ is
quasiconvex in $H$. 

\item \cite{N2} If $\phi\in Aut(G)$ and $A$ is a subgroup of $G$,
 then  $A$ is quasiconvex in $G$ 
if and only if $\phi(A)$ is quasiconvex in $G$.

\item \cite{Swa} Let $(L,X)$ be an automatic structure on $G$ and let $H$ be a
subgroup of $G$. Then $H$ is quasiconvex in $G$ if and only if $H$ is
rational with respect to $L$, that is the pre-image $L_H$ of $H$ in
$L$ is a regular language.

\item  Let $H$ be a quasiconvex subgroup of $G$. Let $X$ be a finite
generating set for $G$ and let $Y$ be a finite generating set for $H$.
For every $y\in Y$ fix a word $v_y$ over $X$ representing $y$. Then
there exist $\lambda>0$ with the following properties.
\smallskip
Suppose $W$ is a $d_Y$-geodesic word over $Y$ representing some
element of $H$. Let $w$ be the $X$-word obtained from $W$ by replacing
each $y\in Y$ by $v_y$. Then $w$ is $(\lambda, \lambda)$-quasigeodesic
with respect to $d_X$. (This follows immediately from the definition
of quasiconvexity). 

\item \cite{K} If $H$ is a quasiconvex subgroup of a word hyperbolic group $G$
then the membership problem for $H$ is solvable.

\endroster
\endproclaim

If $(X,d_X)$, $(Y, d_Y)$ are metric spaces, we say that a map
$f:X\longrightarrow Y$ is a {\it quasi-isometry} if there is $C>0$
such that
\roster 
\item for every $x, x'\in X$ we have
$$\frac{1}{C}d_Y(f(x), f(x')) -C \le d_X(x, x') \le C\cdot d_Y(f(x), f(x'))+C;$$
\item for any $y\in Y$ there is $x\in X$ such that $d_Y(y, f(x))\le C$
(that is the $C$-neighborhood of $f(X)$ is equal to $Y$).
\endroster
\smallskip 
More generally, a map $f:(X,d_X)\longrightarrow (Y, d_Y)$ is called a
{\it quasi-isometric embedding} if there is $C>0$ such that condition
(1) above is satisfied. (Note that $f$ need not be injective).
\medskip
Then the definition of quasiconvexity and the fact that the word
metric is left-invariant immediately imply the following.
\proclaim{Proposition 3.3} Let $G$ be a word hyperbolic group and let
$H$ be a subgroup of $G$. Then the following conditions are
equivalent.
\roster
\item The subgroup $H$ is quasiconvex in $G$.
\item The subgroup $H$ is finitely generated and for some finite
generating set $Y$ of $H$ and for some finite generating set $X$ of
$G$ the inclusion
$$i:(H,d_Y)\longrightarrow (G,d_X)$$
is a quasi-isometric embedding.
\item The subgroup $H$ is finitely generated and for any finite
generating set $Y$ of $H$ and for any finite generating set $X$ of
$G$ the inclusion
$$i:(H,d_Y)\longrightarrow (G,d_X)$$
is a quasi-isometric embedding.
\endroster
\endproclaim
  
Finally, we would like to point out some of the more interesting
properties of quasiconvex subgroups.

\proclaim{Proposition 3.4} Let $G$ be a word hyperbolic group.
\roster
\item \cite{KS, Theorem 1(3)}, \cite{MT} If $H$ is an infinite quasiconvex subgroup of $G$ then $H$ has
finite index in its virtual normalizer

$$V_G(H)=\{g\in G\mid |H:H\cap gHg^{-1}|<\infty, |gHg^{-1}:H\cap
gHg^{-1}|<\infty \}.$$
\item \cite{KS, Theorem 1(2)} If $H$ is an infinite quasiconvex subgroup of $G$ then there are
only finitely many subgroups $L\le G$ such that $L$ contains $H$ as a
subgroup of finite index.
\item \cite{KS, Theorem 1(1)} If $H_1$ and $H_2$ are infinite quasiconvex subgroups of $G$
such that $H=H_1\cap H_2$ has finite index in both $H_1$ and $H_2$
then $H$ has finite index in the subgroup generated by the union of
$H_1$ and $H_2$.
\item \cite{KS, Theorem 3}, \cite{MT} Suppose $H$ is a quasiconvex subgroup of $G$. Suppose $L$ is an
infinite subgroup of $H$ such that $L$ is normal in some subgroup $K$
of $G$ which contains $H$, that is
$$L\le H\le K\le G,$$
$$L \triangleleft K.$$
Then $L$ has finite index in $H$ and $H$ has finite index in $K$.
\item \cite{MT} If $H$ is an infinite quasiconvex subgroup of a word
hyperbolic group $G$ and $t^{-1}Ht\subseteq H$ for some $t\in G$ then $t^{-1}Ht= H$.
\endroster
\endproclaim

\head  4.The combination theorem. \endhead

Suppose $G=\pi_1({\Bbb A}, T)$ where ${\Bbb A}$ is a graph of groups
(see \cite{Ba} and \cite{Ser} for definition and properties of graphs
of groups) with underlying finite graph $A$, maximal subtree $T$ and such that
every vertex group is word hyperbolic and every edge monomorphism is a
quasi-isometric embedding (that is the image of an edge group under the edge
monomorphism is quasiconvex in the appropriate vertex group). 

\definition{Question} When is $G$ word hyperbolic?
\enddefinition
Partial results in this direction can be found in \cite{Gr},
\cite{BGSS}, \cite{KM}, \cite{Gi} and some
other papers. However, the most complete answer to this question up to
date is
given in the Combination Theorem of M.Bestvina and M.Feighn
\cite{BF1}, \cite{BF2}. Before
formulating their result we need to introduce some definitions.

\definition{Definition 4.1 (Annulus)} Let ${\Bbb A}$ be a graph of groups with vertex groups $A_v$, $v\in VA$ and edge groups $A_e$, $e\in EA$. For an edge $e\in EA$ denote the edge homomorphisms by $\alpha_e : A_e\longrightarrow A_{\partial_0(e)}$ and $\alpha_{\overline e} : A_e\longrightarrow A_{\partial_1(e)}$.
\smallskip
A {\it combinatorial annulus} of length $2M+1$ is a diagram $\Sigma$
as in Figure 4.1

\topinsert
\epsfxsize=15cm\epsfbox{f4.1.eps}
\captionwidth{7cm}
\botcaption{Figure 4.1 (Annulus)}
\endcaption
\endinsert

satisfying the following properties:
\roster
\item the sequence $(e_{-M},e_{-M+1},\dots ,e_0,\dots ,e_M)$ is an
edge-path in $A$;
\item for every $i=-M,\dots M-1$ we have $$a_i\in A_{v_i}$$ where
$v_i={\partial_1(e_i)}={\partial_0(e_{i+1})}$;
\item for every $i=-M, \dots , M$ we have
$$c_i\in A_{e_i}, c_i\ne 1$$ and

$$ {a_i}^{-1} \alpha_{\overline {e_i}}(c_i) a_i =\alpha_{e_{i+1}}(c_{i+1})$$
\endroster

Thus in $G$

$$\align &(e_{-M}a_{-M}e_{-M+1}\dots a_{M-1}e_M)^{-1}\cdot \alpha_{e_{-M}}(c_{-M})\cdot e_{-M}a_{-M}e_{-M+1}\dots a_{M-1}e_M= \\
& \alpha_{{\overline e}_M}(c_M)
\endalign$$
and 
$$\align &\alpha_{e_{-M}}(c_{-M})\cdot e_{-M}a_{-M}e_{-M+1}\dots a_{M-1}e_M= \\
&=e_{-M}a_{-M}e_{-M+1}\dots a_{M-1}e_M\cdot \alpha_{{\overline e}_M}(c_M)
\endalign$$

\enddefinition

\definition{Definition 4.2 (Essential annulus)}
An annulus $\Sigma$ as above is called {\it essential} if the sequence
$$(e_{-M},a_{-M},e_{-M+1},\dots ,a_{M-1},e_M)$$ contains no "pinches"
that is it has no subsequences of the form
$$ e, \alpha_{\overline e} (s), \overline e$$
where $s\in A_e$, and $e$ is an edge of $A$.
\enddefinition
\remark{Remark 4.3} From the geometric viewpoint if $\Sigma$ is an
essential annulus then the sequence $(e_{-M},a_{-M},e_{-M+1},\dots
,a_{M-1},e_M)$ represents an edge-path without backtracks of length
$2M+1$ in the Bass-Serre universal covering tree $\hat T$ corresponding to the
graph of groups ${\Bbb A}$. Moreover, the definition of an annulus in
this case is just an algebraic restatement of the fact that there is a
nontrivial element in the fundamental group $\pi_1({\Bbb A}, \ast)$ of
the graph of groups ${\Bbb A}$ which fixes this edge-path pointwise.
\endremark

\smallskip
\definition{Definition 4.4} Let ${\Bbb A}$ be a graph of groups as above.
Fix generating sets for all the vertex groups $A_v$, all the
edge groups $A_e$ and the word metrics $d_v$, $d_e$ induced by them.
Then if $\Sigma$ is an annulus as in Definition 4.1, we say that

\roster
\item the {\it girth} of ${\Sigma}$ is $d_{e_0}(c_0)$;
\item the {\it width} of $\Sigma$ is $\max \{d_{v_i}(a_i) | -M \le i<M\}$.
\endroster

Let $\lambda >1$. An annulus $\Sigma$ is called $\lambda$-{\it
hyperbolic} (see Figure 4.2) if

$$\lambda d_{e_0}(c_0)\le \max\{ d_{e_M}(c_M), d_{e_{-M}}(c_{-M})\}.$$ 
\enddefinition

\topinsert
\epsfxsize=8cm\epsfbox{f4.2.eps}
\captionwidth{9cm}
\botcaption{Figure 4.2 (Hyperbolic annuli)}
\endcaption
\endinsert

\proclaim {Theorem 4.5 (Combination Theorem)} (M.Bestvina and M.Feighn \cite{BF2})
Suppose $G=\pi_1({\Bbb A}, T)$ where ${\Bbb A}$ is a graph of groups
with underlying finite graph $A$, maximal subtree $T$ and such that
every vertex group is word hyperbolic. Suppose also that every edge monomorphism
$\alpha_e: A_e\longrightarrow A_{\partial_0(e)}$ is a
quasi-isometric embedding (that is $\alpha_e(A_e)$ is quasiconvex in $A_{\partial_0(e)}$).
Fix finite generating sets for all $A_v$, $A_e$ and the word metrics $d_v$, $d_e$ induced by them.
\smallskip

Suppose there exist $\lambda >1$ and $M \ge 1$ such that the following holds.
For each $\rho >0$ there is $H(\rho)$ such that any essential annulus
of length $2M+1$, width at most $\rho$ and girth at least $H(\rho)$ is $\lambda$-hyperbolic.
\medskip
Then $G$ is word hyperbolic.
\endproclaim
\remark{Remark 4.6} Suppose that the conditions of the Combination Theorem
as stated above are satisfied for a particular choice of generating
sets for the vertex groups $A_v$ and edge groups $A_e$. Then it is obvious from the definition
of an annulus that for any integer $k>0$ every essential annulus of
width at most $\rho$ and girth at least $H(\rho)$ and of length
$2kM+1$ is $\lambda^k$-hyperbolic. Therefore the conditions of the
Combination Theorem for the graph of groups ${\Bbb A}$ are satisfied
(with different $M$ and $\lambda$) for any other choice of finite
generating sets for the vertex groups and the edge groups of ${\Bbb
A}$.
\endremark

There is one easy corollary of the Combination Theorem that is
extremely useful when working with hyperbolic groups.

\proclaim{Corollary 4.7} Let ${\Bbb A}$ be a finite graph of groups such
that all vertex groups are word hyperbolic and the images of the edge
groups under edge monomorphisms are quasiconvex in the appropriate
vertex groups. Let $T$ be a maximal tree in the underlying graph $A$ of
${\Bbb A}$ and let $G=\pi_1({\Bbb A}, T)$ be the fundamental group of
the graph of groups ${\Bbb A}$ with respect to $T$. 

\smallskip 
Suppose there is an integer $M>0$ such that there are no essential
annuli of length $2M+1$ corresponding to ${\Bbb A}$.
\medskip
Then $G$ is word hyperbolic.
\endproclaim
A direct restatement of Corollary 4.7 yields the following.

\proclaim{Corollary 4.8} Let ${\Bbb A}$ be a finite graph of groups such
that all vertex groups are word hyperbolic and the images of the edge
groups under edge monomorphisms are quasiconvex in the appropriate
vertex groups. Let $T$ be a maximal tree in the underlying graph $A$ of
${\Bbb A}$ and let $G=\pi_1({\Bbb A}, T)$ be the fundamental group of
the graph of groups ${\Bbb A}$ with respect to $T$. Let $\hat T$ be
the Bass-Serre universal covering tree of ${\Bbb A}$ on which $G$
acts.
\smallskip 
Suppose there is an integer $M>0$ such that no nontrivial element of
$G$ fixes a segment of length $M$ in $\hat T$ pointwise. 
\medskip
Then $G$ is word hyperbolic.
\endproclaim

We would like to point out two explicit and most often used
applications of the above results. First recall the following definition.
\definition{Definition 4.9} A subgroup $H$ of a group $G$ is called {\it
malnormal} if for every $g\in G$, $g\not\in H$
$$g H g^{-1} \cap H=1$$
\enddefinition

\proclaim{Proposition 4.10} Let $A$ and $B$ be word hyperbolic groups such
that the subgroup $C=A\cap B$ is quasiconvex in both $A$ and $B$.
Suppose also that $C$ is malnormal in $B$. Then the amalgamated free
product $G=A\ast_C B$ is word hyperbolic.
\endproclaim
\demo{Proof}
Let $\Bbb Y$ be the edge of groups with the vertex groups $A$, $B$ and
the edge group $C$ as shown in Figure 4.3.
\topinsert
\epsfxsize=7cm\epsfbox{f4.3.eps}
\captionwidth{5cm}
\botcaption{Figure 4.3}
\endcaption
\endinsert

Thus $Y$ has the edge $e$ with the endpoints $\partial_0(e)=v_0$,
$\partial_1(e)=v_1$. The vertex groups for $\Bbb Y$ are $Y_{v_0}=A$,
$Y_{v_1}=B$ and the edge groups are $Y_e=Y_{\overline e}=C$. The
boundary monomorphisms $\alpha_e:C\rightarrow A$ and
$\alpha_{\overline e}:C\rightarrow B$ are the inclusions. Then
obviously $G=A\ast_C B\cong \pi_1({\Bbb Y}, v_0)$.
\medskip
We claim that there are no essential annuli of length 3 for the graph
of groups  $\Bbb Y$. Indeed, suppose that $\Sigma$ is an essential
annulus of length 3. Then up to axial symmetry $\Sigma$ has the form
as in Figure 4.4, (that is the sequence of edges in the top label
of $\Sigma$ is $\overline e, e, \overline e$ rather than $e, \overline
e, e$).
\midinsert
\epsfxsize=8cm\epsfbox{f4.4.eps}
\captionwidth{5cm}
\botcaption{Figure 4.4}
\endcaption
\endinsert

Then $c_{-1}, c_0, c_1\in C-\{1\}$, $a_{-1}\in A$, $b_0\in B$. Since the
boundary monomorphisms are just inclusions, by the definition of an
annulus we have

$$b_0^{-1}c_0b_0=c_1\text { in } B.$$
Since $C$ is malnormal in $B$ and $c_0\ne 1$, this implies $b_0\in C$.
Therefore the sequence $e, b_0, \overline e$ is a ``pinch'' and the
annulus $\Sigma$ is not essential. This contradicts our assumptions.
\smallskip

Thus there are indeed no essential annuli of length 3 and therefore
the group $G$ is word hyperbolic by Corollary 4.7.

\enddemo

\proclaim{Proposition 4.11}

Let $G$ be word hyperbolic group. Let $C_1$ and $C_2$ be isomorphic
quasiconvex subgroups of $G$ and let $\phi:C_1\rightarrow C_2$ be an
isomorphism. Suppose that both $C_1$ and $C_2$ are malnormal on $G$
and, moreover, no nontrivial element of $C_1$ is conjugate in $G$ to
an element of $C_2$.
\smallskip

Then the HNN-extension 

$$K=<G,t | t^{-1}c t=\phi(c), c\in C_1>$$

is word hyperbolic.
 
\endproclaim
\demo{Proof}
Once again, it is easy to see that if ${\Bbb A}$ is the loop of groups
associated to the HNN presentation of $K$ then there are no essential
annuli of length 3 corresponding to $\Bbb A$ (the details are left to
the reader). Therefore by Corollary 4.7 the group $K$ is word hyperbolic.
\enddemo
\remark{Remark 4.12} Suppose that $\Bbb A$ is a finite graph of groups
such that all the edge groups are infinite cyclic. Let $T$ be a
maximal tree in $A$ and let $G=\pi_1({\Bbb A}, T)$. It is not hard to
see that if $\Bbb A$ has essential annuli of an arbitrary big
length, then $G$ contains a so-called Baumslag-Solitar subgroup, that
is a subgroup of the form
$$B(m,n)=<a,t| t^{-1}a^nt=a^m>$$ for some $m\ne 0, n\ne
0$. Baumslag-Solitar groups can never be subgroups of word hyperbolic
groups (see Proposition 4.4(6)).
Therefore if ${\Bbb A}$ is a finite graph of groups with word
hyperbolic vertex groups and infinite cyclic edge groups then the
fundamental group $G$ of the graph of groups ${\Bbb A}$ is word
hyperbolic if and only if $G$ does not contain Baumslag-Solitar
subgroups \cite{BF2}.
\endremark 
We will later need the following simple lemma
\proclaim{Lemma 4.13} Let $\Bbb A$ be a graph of groups with
and underlying graph $A$. Suppose $e$ is an edge of $A$ such that the
subgroup $\alpha_e(A_e)$ is malnormal in the vertex group $A_v$ where
$v=\partial_0(e)$. Then there are no essential annuli $\Sigma$
corresponding to $\Bbb A$ such that the top label of $\Sigma$ contains
a subsequence $\overline e, a, e$ where $a\in A_v$.
\endproclaim
\demo{Proof}
Suppose, on the contrary, such an essential annulus $\Sigma$ exists.
Then $\Sigma$ contains a sub-annulus of the form shown in Figure 4.5.

\topinsert
\epsfxsize=7cm\epsfbox{f4.5.eps}
\captionwidth{5cm}
\botcaption{Figure 4.5}
\endcaption
\endinsert

Here $a\in A_v$, $c,c'\in A_e=A_{\overline e}$, $c\ne 1, c'\ne 1$.
Then by the definition of an annulus
$$a^{-1}\alpha_e(c) a =\alpha_e(c').$$
Since  $\alpha_e(A_e)$ is malnormal in $A_v$, this implies that $a\in
\alpha_e(A_e)$. However, this contradicts our assumption that the
annulus $\Sigma$ is essential.
\enddemo

\head 5.Proofs of the main results. \endhead

\proclaim{Theorem 5.1 (Theorem A)} Let $G$ be a non-elementary
torsion-free word hyperbolic group. Then there exists another word
hyperbolic group $G^*$ such that $G$ is a subgroup of $G^*$ and $G$ is
not quasiconvex in $G^*$.
\endproclaim

\demo {Proof of Theorem A}

Let $G$ be as in Theorem 5.1. Then by Theorem 7.7 (which is proved
later in section 7) there exists a subgroup $F$ of $G$
such that
\roster 
\item $F$ is free of rank 2
\item $F$ is quasiconvex in $G$
\item $F$ is malnormal in $G$.
\endroster
\medskip
Say $F=F(a,b)$ is free on $a,b\in G$. For any $f\in F$ we will denote
by $l_F(f)$ the length of the freely reduced word in $a,b$ representing
$f$.
\medskip
Take $\phi: F(a,b)\longrightarrow F(a,b)$ to be an endomorphism such that
\roster
\item "(a)" the subgroup $im(\phi)=\phi(F(a,b))$ is malnormal in $F(a,b)$ (and so in $G$);
\item "(b)" the map $\phi$ is {\it uniformly length-expanding} that is
for every $f\in F(a,b)$
$$l_F(\phi(f)) \ge 2 l_F(f).$$
(Note that property (b) implies that $\phi$ is a monomorphism). 
\endroster
\medskip
Such a $\phi$ always exists. 
For example, we can take $\phi(a)=a b^3 a$ and $\phi(b)=b a^3 b$. Note that the word $ab^3a$ begins with $a$ and ends with $a$ and
the word $b a^3 b$ begins with $b$ and ends with $b$. Therefore  
for each $f\in F(a,b)$ we have $l_F(\phi(x)) =5 l_F(f)$. Also the subgroup
$gp(a b^3 a, b a^3 b)$ is malnormal in $F(a,b)$ by Lemma 6.4.
\medskip
So assume that $\phi$ is any endomorphism of $F(a,b)$ that satisfies
conditions (a), (b) above.
\smallskip
Put $$G^*=<G,t | t^{-1} a t=\phi (a), t^{-1} b t =\phi (b)>$$ to be the HNN-extension of $G$ along $\phi$.
\smallskip
{\bf Claim} Then
\roster
\item $G^*$ is word hyperbolic;
\item $G$ is not quasiconvex in $G^*$.
\endroster

To see that $G^*$ is word hyperbolic, note that $G$ is isomorphic to
the fundamental group of the graph $\Bbb A$ of groups shown in Figure 4.6.

\topinsert
\epsfxsize=8cm\epsfbox{f4.6.eps}
\captionwidth{5cm}
\botcaption{Figure 4.6}
\endcaption
\endinsert

More precisely, the underlying graph $A$ of $\Bbb A$ has a single
vertex $v$ and an edge $e$ with $\partial_0(e)=\partial_1(e)=v$.
Also the vertex and edge groups are $A_v=G$, $A_e=F(x,y)$ where
$F(x,y)$ is a free group of rank 2 with basis $x,y$. The boundary
monomorphisms for $\Bbb A$ are defined as follows:

$$\alpha_e(x)=a, \alpha_e(y)=b$$
 and
$$\alpha_{\overline e}(x)=\phi(a), \alpha_{\overline e}(y)=\phi(b).$$

Then obviously $G^*\cong \pi_1({\Bbb A}, v)$.
We will use the combination theorem to show that $G^*$ is word hyperbolic.
Note first that every essential annulus $\Sigma$ is {\it one-directed}, that is the label of its upper boundary does not contain subsequences of the type

$$ e, a, \overline e \quad \text{or}\quad \overline e, a, e$$
This follows by Lemma 4.13 because subgroups $F(a,b)$ and
$\phi(F(a,b))$ are malnormal in $G$ and the annulus $\Sigma$ is essential.
\smallskip
Fix a finite generating set $\Cal G$ for $G$ and the word metric
$d_{\Cal G}$ corresponding to $\Cal G$. Recall that for any $f\in
F(a,b)$ we denote the freely reduced length of $f$ in $a, b$ by
$l_F(f)$. Also if $c\in F(x,y)=A_e$, then we will denote the freely
reduced length of $c$ in $x,y$ by $l_X(c)$.
\smallskip 
Put $\lambda =\frac {3}{2}$ and $M=1$.
Now for every $\rho >0$ let $q(\rho )$ be the maximum of $l_F$-lengths
of those elements of $F(a,b)$ whose $\Cal G$-length is at most $\rho$.
Put $H(\rho)=20 q(\rho )$.
\smallskip
Let $\rho>0$ be an arbitrary number.
We claim that for the graph of groups $\Bbb A$ every essential annulus
$\Sigma$  of length 3, width at most $\rho$ and girth at least $H(\rho)$ is $\lambda$-hyperbolic.
\smallskip
Indeed, since $\Sigma$ is one-directed, up to an axial symmetry we can
assume that 
$\Sigma$ looks like the annulus shown in Figure 4.7. (That is the upper label of $\Sigma$
is of the form $e,a_{-1}, e, a_0, e$.) 
\topinsert
\epsfxsize=8cm\epsfbox{f4.7.eps}
\captionwidth{5cm}
\botcaption{Figure 4.7}
\endcaption
\endinsert

Here $a_{-1}, a\in G$, $c_{-1}, c_0, c\in F(x,y)-\{1\}$ and 
$$a_i^{-1}\alpha_{\overline e}(c_i)a_i=\alpha_e(c_{i+1})$$
for $i=-1,0$.
Since $F(a,b)$ is malnormal in $G$ and the subgroups $\alpha_e(A_e)$
and $\alpha_{\overline e}(A_e)$ are contained in $F(a,b)$, this
implies that $a_0,a_{-1}\in F(a,b)$.
Also $l_{\Cal G}(a_{-1}), l_{\Cal G}(a_0)\le \rho$ and therefore by
the choice of $q(\rho)$ we have  $l_F(a_i)\le q(\rho )$, $i=-1,0$.

Note also that 

$l_F(\alpha_{\overline{e}}(c_0))\ge 2 l_X(c_0)$,  
$l_F(\alpha_{e}(c_1))=l_X(c_1)$ and $l_X(c_0)\ge H(\rho)=20q(\rho )$.

Thus $$\align &l_X(c_1)=l_F(\alpha_{\overline e}(c_1))\ge \\
&l_F(\alpha_{\overline e}(c_0)) -2q(\rho) \ge 2 l_X(c_0)-2q(\rho )\ge \frac{3}{2} l_X(c_0).\endalign$$

Thus $\Sigma$ is $\dsize \frac{3}{2}$-hyperbolic and the group $G^*$ is word hyperbolic by the combination theorem.
\medskip
We will now show that $G$ is not quasiconvex in $G^*$.
Assume that, on the contrary, $G$ is quasiconvex in $G^*$. Since
$F(a,b)$ is quasiconvex in $G$ this implies that $F(a,b)$ is
quasiconvex in $G^*$ (see Proposition 3.2(6)). Denote by $d_*$ the word metric on $G^*$
corresponding to the generating set ${\Cal C}\cup \{t\}$.
Then there is $C>0$ such that for every $f\in F(a,b)$

$$l_F(f)\le C l_*(f)$$

Take $f_n=t^{-n} a t^n=\phi^n(a)\in F(a,b)$, $n>0$.
Then $l_F(f_n) \ge 2^n$ by the properties of $\phi$ and $l_*(f_n)\le 2n+l_*(a)$.
Therefore

$$2^n \le l_F(f_n) \le C l_*(f_n) \le 2C\cdot n +Cl_*(a)$$
for every integer $n >0$.
This gives us a contradiction.
Therefore $G$ is not quasiconvex in $G^*$ and Theorem A is proved.
\enddemo
The remainder of this paper is devoted to proving Theorem 7.7 which
was used in the proof of Theorem A and which asserts that a
non-elementary torsion-free word hyperbolic group always has a free
quasiconvex malnormal subgroup.

\head 6.Malnormality in free groups. \endhead

In order to prove the existence of a free malnormal quasiconvex
subgroup in a non-elementary torsion-free hyperbolic group we need to
accumulate a certain amount of information about malnormal subgroups
of free groups. 

Our main goal in this section is to prove the following theorem (the
problem of finding a malnormal subgroup in a hyperbolic group will be
reduced to this statement).

\proclaim{Theorem C}
Let $F=F(X)$ be a nonabelian free group of finite
rank with basis $X$.  Let
$H_1, \dots H_k$ be finitely generated subgroups of $F$ of infinite
index. Then there exist elements $a,b\in F$ such that
\roster
\item the subgroup $H=gp(a,b)$ is a free group of rank two which is
malnormal in $F$;
\item no nontrivial element of $H$ is conjugate in $F$ to an element
of $H_1\cup\dots \cup H_k$.
\endroster
\endproclaim
In order to prove Theorem C we will need the following series of lemmas.
\proclaim{Lemma 6.1} Let $F=F(X)$ be a non-abelian finitely generated free group
with a finite basis $X$. Let $H_1,\dots , H_k$ be finitely generated
subgroups of $F$ of infinite index.

Then there exists an element $y\in F$ such that no nontrivial power of $y$ is
conjugate to an element of $H_1\cup \dots \cup H_k$.
\endproclaim

\demo{Proof}

For each subgroup $H_i$ we choose a Nielsen reduced basis $B_i$.
Let $K$ be the maximal length of the elements of $B_1\cup \dots \cup
B_k$. Note that $K$ obviously has the following property. If $w$ is a
freely reduced word in $X$ representing an element of $H_i$ and $v$ is
a subword of $w$ then there are words $a$, $b$ in $X$ of length at
most $K$ such that $avb\in H_i$.

We now claim that there is a nontrivial freely reduced word $y$ in $X$
which is not a subword of any freely reduced word $w$ representing an
element of $H_1\cup H_2 \cup \dots \cup H_k$.

Indeed, if there are no such $y$ then $F=F(X)$ is the following finite union: 

$$F(X)=\bigcup_{i=1}^{i=k}\ \bigcup_{l(a), l(b)\le K} a H_i b.$$

Note that every set $a H_i b$ is in fact a right coset $(a H_i a^{-1})(ab)$.
However all the subgroups $H_1, \dots H_k$ (and so all of their
conjugates) are of infinite index in
$F$. Thus $F$ is covered by finitely many right cosets of subgroups of
infinite index. By the result of B.H.Neumann \cite{N1},
this is impossible for any group which gives us a contradiction.
Thus the claim is proved and there exists a freely reduced $y$ which
is never a subword of a freely reduced word representing an element of
$H_i$. 

Moreover, we may assume that $y$ is cyclically reduced. Suppose, on
the contrary, that $y$ is not
cyclically reduced and has the form $y=x^{\epsilon}y' x^{-\epsilon}$
where $\epsilon =\pm 1$ and $x$ is an element of $X$. Then we take any
other letter $z\in X$, $z\ne x$ (it is possible to do since we assumed
$F=F(X)$ is non-abelian and so has rank at least two). Now replace $y$
by $y_1=zyz$. 
Obviously, $y_1$ is cyclically reduced and $y_1$ is not a subword of a
freely reduced word representing an element of $H_i$. 

Thus we may indeed assume that $y$ is cyclically reduced to begin with.
We now claim that no nontrivial power of $y$ is conjugate to an
element of $H_i$. Indeed, suppose $q^{-1}y^p q =h\in H_i$ for some
$p> 0$, $q\in F$, $1\le i\le k$. 

Take $m>0$ to be such that $l(q)<l(y^{pm})=pml(y)$.
Then $q^{-1}y^{3pm} q= h^{3m}\in H_i$. It follows from the choice of
$m$ that the freely reduced word representing $q^{-1}y^{3pm} q= h^{3m}$
contains $y$ as a subword which contradicts the properties of $y$.

Lemma 6.1 is proved.
\enddemo

\proclaim {Lemma 6.2} Let $F=F(X)$ be a free nonabelian group of finite
rank with basis $X$. Let $a$ and $b$ be nontrivial cyclically reduced
words in $X$ such that no nontrivial power of $a$ is conjugate to a
power of $b$. Then there exists $K>0$ such that for any $n\ne 0, m\ne
0$ the length of the maximal subword of $a^n b^m$, that
freely reduces to the identity, is at most $K$.
\endproclaim
\demo{Proof}
Let $N$ be the number of distinct elements of $F$ of length at most
$l(b)$. Put $K=Nl(a)=l(a^N)$. Suppose that in some product $a^n b^m$ a
terminal segment $q$ of $a^n$ of length at least $K$ is freely cancelled with
an initial segment of $b^m$.
Then $q=a^{\pm N}$. This means that $a^{\pm N}$ is an initial segment
of $b^n$. We may assume that in fact $a^N$ is an initial segment of
$b^n$. Put $e_0=1$. For every $i=1,\dots , N$ there is an element
$e_i\in F$ of length at most $l(b)$ such that $a^ie_i$ is a power of
$b$. All $N+1$ elements $e_0,e_1,\dots, e_N$ have length at most
$l(b)$. Therefore by the choice of $N$ there are $i_1, i_2$, $0\le
i_1<i_2\le N$ such that $e_{i_1}=e_{i_2}=e$. Thus $a^{i_1}e=b^{n_1}$
and $a^{i_2}e=b^{n_2}$ for some $n_1,n_2$. Therefore
$e^{-1}a^{i_2-i_1}e=b^{n_2-n_1}$. However, this contradicts our
assumption that no power of $a$ is conjugate to a power of $b$.
Thus we have proved that the length of the maximal subword of $a^n b^m$, that
freely reduces to the identity, is at most $2K=2Nl(a)$.
Lemma 6.2 is proved.
\enddemo

\proclaim{Lemma 6.3}Let $F=F(X)$ be a non-abelian finitely generated free group
with a finite basis $X$. Let $H_1,\dots , H_k$ be finitely generated
subgroups of $F$ of infinite index.
Let $y\in F$ be a cyclically reduced word in $X$ such that no
nontrivial power of $y$ is
conjugate to an element of $H_1\cup \dots \cup H_k$.
Let $t$ be another cyclically reduced word in $X$ such that no
nontrivial power of $t$ is conjugate to a power of $y$.

Then there exists integer $m>0$ with the following property. If $n\ge
m$ then no nontrivial element of the subgroup $D_n=gp(y^n t^{3n} y^n,
t^n y^{3n} t^n)$ is conjugate to an element of $H_1\cup \dots \cup H_k$.

\endproclaim
\demo{Proof of Lemma 6.3}

Let $K_1>l(y)+l(t)$ be such that for any $k\ne 0, s\ne 0$ the length
of the maximal subword of $y^kt^s$ that freely reduces to the identity
is at most $K_1$. The existence of such a constant $K_1$ follows from
Lemma 6.2.
For each subgroup $H_i$ we choose a Nielsen reduced basis $B_i$.
Let $K_2$ be the maximal length of the elements of $B_1\cup \dots \cup
B_k$. Put $K=max(K_1, K_2)$. 
Note that by the choice of $y,t$ the subgroup $D=gp(y,t)\le F$ is free
of rank two and therefore for every $n>0$ the group $D_n=gp(y^n t^{3n} y^n,
t^n y^{3n} t^n)$ is also a free group of rank two with basis $y^n t^{3n} y^n,
t^n y^{3n} t^n$.
\smallskip
Let $m_1>0$ be such that $l(y^{m_1})=m_1 l(y)> K \ge K_1$ and
$l(t^{m_1})=m_1l(t)>K\ge K_1$. Then it is clear
that for every $n\ge 3m_1$ every freely reduced word $w$ representing
a nontrivial element of $D_n$ contains $y^{n-2m_1}$ as a subword.
Moreover, for every $n\ge 3m_1$, any $z\in F$ and any nontrivial
element $d\in D_n$ some sufficiently high power of $z^{-1}dz$, when written
in the freely reduced form, contains $y^{n-2m_1}$ as a subword. 
\smallskip
Let $N$ be the number of distinct elements of $F$ of length at most $K$.  
Put $$m=3m_1+N+1.$$ Suppose now that $n\ge m$ and $d\in D_n, d\ne 1$.
\smallskip
We claim that no nontrivial power of $d$ is conjugate in $F$ to an
element of $H_i$. Indeed, suppose $z^{-1} d^p z=h\in H_i$ for some
$z\in F$, $p>0$. Then by the remarks above some sufficiently high
power $h^s$ of $h$ contains, when written in a freely reduced form, a
subword $y^{n-2m_1}$.
\medskip
Note also that by the properties of a Nielsen basis for $H_i$, if $v$
is an initial segment of a freely reduced word in $X$ representing an
element of $H_i$, then for some element $a\in F$ with $l(a)\le K_2\le
K$ we have $va\in H_i$.
\medskip
Let the freely reduced word $w$ representing $h^s$ have the form
$w\equiv v y^{n-2m_1} v'$.  Note that $n\ge m=3m_1+N+1$ and therefore
$n-2m_1\ge N+1$. Thus there are some elements $a_0, a_1,\dots a_N\in
F$ with $l(a_j)\le K_2\le K$ such that $vy^ja_j\in H_i$ for every
$j=0,1,\dots N$. Since there are only $N$ different elements of length
at most $K$ in $F$, there exist $j_1<j_2$, $0\le j_1, j_2\le N$
such that $a_{j_1}=a_{j_2}=a$.
Thus we have $vy^{j_1}a\in H_i$ and $vy^{j_2}a\in H_i$ which implies
$a^{-1}y^{j_2-j_1}a\in H_i$. This contradicts our choice of $y$.
\smallskip
We have shown that if $n\ge m=3m_1+N+1$ then no nontrivial element of
$D_n$ is conjugate to an element of $H_1\cup \dots \cup H_k$.
Lemma 6.3 is proved.

\enddemo

\proclaim {Lemma 6.4} Let $F(y,t)$ be the free group on $y, t$. Then
the subgroup $$M=gp(y t^3 y, t y^3 t)$$ is malnormal in $F(y,t)$.
\endproclaim
\demo{Proof}
This fact can be established by elementary means, however the proof is
a rather lengthy exercise. We have verified the correctness of Lemma 6.4
using the computational group theory software package MAGNUS \cite{Ma}.
\enddemo

\proclaim{Lemma 6.5} Let $F=F(a,b)$ be the free group on $a,b$ and let 
$n>0$, $m>0$
be some integers. Let $L=gp(a^n, b^m)$ be the subgroup of $F$ generated
by $a^n, b^m$. Suppose that $z h z^{-1} =h'$ for some $h, h'\in L$,
$h\ne 1$, $z\in F-L$. Then $h$ is conjugate in $L$ either to a power
of $a^n$ or to a power of $b^m$.
\endproclaim
\demo{Proof}
First notice that it is enough to prove Lemma 6.5 for those $z$ which
are shortest in their coset classes $zL$. Indeed, suppose that Lemma 6.5
has been established for such $z$.
Now let $z$ be an arbitrary element as in Lemma 6.5. Then for some
$h_0\in L$ we have $z=z_1h_0$ where $z_1$ is shortest in the coset
class $zL$.
Then $z h z^{-1} =h'$ implies that $z_1 (h_0 h h_0^{-1}) z_1^{-1}
=h'$. Therefore by our assumption $h_0 h h_0^{-1}$ is conjugate in $L$
to a power of $a^n$ or a power of $b^n$ and the same is true for $h$.
\medskip
Thus from now on we will assume that $z$ is shortest in $zL$.
Without a loss of generality we may assume that the last letter of $z$
(when written as a freely reduced word in $a,b$) is $a^{\pm 1}$.
If $z$ is a power of $a$ then the statement of Lemma 6.5 is obvious. If
$z$ is not a power of $a$ then $z\equiv qb^k a^s$ where $k\ne 0$,
$s\ne 0$. Moreover, since $z$ is shortest in $zL$, it cannot be
shortened by multiplying on the right by $a^{\pm n}$. Therefore $0<
|s|<n$.
\medskip
If $h$ is a power of $a^n$ or a power of $b^m$, then there is nothing
to prove. Suppose therefore that $h$ is not a power of $a^n$ and that
$h$ is not a power of $b^m$. Then $h=x_1 x_2\dots x_p$, where $p\ge
2$, each of $x_i$ is a nonzero power of $a^n$ or $b^m$ and the
sequence strictly alternates.
\medskip
There are several cases to consider.
\smallskip
{\bf Case 1} Both $x_1$ and $x_p$ are powers of $b^m$. Then $z
hz^{-1}\equiv qb^k a^s x_1 x_2 \dots x_p a^{-s} b^{-k} q^{-1}$ is
freely reduced as written. Since $0<|s|<n$, this implies
$zhz^{-1}\not\in L=gp(a^n,b^m)$ which contradicts our assumptions.
\smallskip
{\bf Case 2} Both $x_1$ and $x_p$ are powers of $a$. Then $p\ge 3$ and
$x_1=a^{n_1n}$, $x_2=a^{n_2n}$. In this case the freely reduced form
of $z hz^{-1}$ is $qb^k a^{n_1n+s} x_2\dots x_{p-1} a^{n_2n-s}b^{-k}
q^{-1}$. Again, since $0<|s|<n$, both numbers $n_1n+s$ and $n_2n-s$
are not divisible by $n$. Therefore $zhz^{-1}\not\in L=gp(a^n,b^m)$
which contradicts our assumptions.
\smallskip
{\bf Case 3} The word $x_1$ is a power of $a^n$ and the word $x_p$ is
a power of $b^m$. Then $x_1=a^{n_1n}$. In this case
the freely reduced form of $z h z^{-1}$ is
$qb^k a^{n_1n+s} \dots x_p a^{-s} b^{-k} q^{-1}$.
Once again, since $0<|s|<n$, we conclude that both $s$ and $n_1n+s$
are not divisible by $n$. Thus $zhz^{-1}\not\in L=gp(a^n,b^m)$
which contradicts our assumptions.
\smallskip
{\bf Case 4} The word $x_1$ is a power of $b^m$ and the word $x_2$ is
a power of $a^n$. This case reduces to Case 3 if we replace the
equality $z hz^{-1}=h'$ by $z h^{-1}z^{-1}=(h')^{-1}$.
\medskip
Thus we proved that $p>1$ is impossible, that is $h$ is a power of
$a^n$ or a power of $b^m$.
This completes the proof of Lemma 6.5.
\enddemo

\proclaim {Lemma 6.6} Let $F=F(a,b)$ be the free group on $a,b$.
Then for any $n>0, m>0$ the subgroup $M=gp(a^n b^{3m} a^n, b^m a^{3n}
b^m)$ is malnormal in $F$.
\endproclaim
\demo{Proof}
Let $L=gp(a^n, b^m)$. Observe first that if $z^{-1} h z= h'$ where $h,
h'\in L-{1}$, $z\in F-L$ then by Lemma 6.5 $h$ is conjugate to either a power of
$a$ or a power of $b$. 
\smallskip
Suppose now that $M$ is not malnormal in $F$. Then for some nontrivial
$h\in M$ and some $z\in F-M$ we have $z^{-1} h z\in M$.
By Lemma 6.4 the subgroup $M$ is malnormal in $L$ and therefore
$z\not\in L$. Then by the previous observation $h$ is conjugate to
either a power of $a$ or a power of $b$. Clearly, neither of these cases are possible
which gives us a contradiction.
\smallskip
Thus $M$ is malnormal in $F$ which completes the proof of Lemma 6.6.
\enddemo

\proclaim{Lemma 6.7} Let $F=F(X)$ be a free group on $X$. Let $v$ and
$s$ be nontrivial freely reduced words in $X$ such that $l(v)>l(s)$.

Suppose also that $sv\equiv v\alpha$, where $l(\alpha)=l(s)$  and the products $sv$ and
$v\alpha$ are freely reduced as written. Let $n$ be the biggest
integer such that $nl(s)\le v$. Then $v\equiv s^n q$ where the
product $s^nq =s s s\dots s q$ is freely reduced as written and
$l(q)< l(s)$.
\endproclaim 
\demo{Proof}
We will prove Lemma 6.7 by induction on $n$.
When $n=1$, the statement is obvious.
Suppose now that $n>1$ and that the statement has been proved for all
smaller values of $n$.
\smallskip
Note that since $sv\equiv v\alpha$ and $l(v)>l(s)$, the word $s$ is a
proper initial segment of $v$. Thus $v$ has the form $v\equiv sv_1$.
Therefore $sv\equiv v\alpha$ implies $ssv_1\equiv sv_1\alpha$.
Hence $s$ is cyclically reduced and $sv_1\equiv v_1\alpha$.
Note that since $n>1$, $l(v)\ge 2l(s)$ and so $l(v_1)\ge l(s)$.
\smallskip
If $l(v_1)=l(s)$ then $v_1=s$, $v=s^2$ and the statement of Lemma 6.7
obviously holds.
\smallskip
If $l(v_1)>l(s)$ then by the inductive hypothesis $v_1\equiv
s^{n-1}q$ with $l(q)<l(s)$. Therefore $v\equiv sv_1\equiv s s^{n-1}q\equiv
s^nq$. Once again the statement of Lemma 6.7 holds.

\enddemo

\proclaim{Lemma 6.8} Let $F=F(X)$ be a free non-Abelian group of finite
rank with basis $X$. Let $a, b$ be cyclically reduced words in $X$
such that no nontrivial power of $a$ is conjugate in $F$ to a power of
$b$. Then there exists some number $M>0$ such that $a^M$ is not a
subword of a power of $b$ and $b^M$ is not a subword of a power of
$a$.
\endproclaim

\demo{Proof}

Let $N$ be the number of elements of $F$ of length at most
$l(b)$. Suppose that $a^N$ is a subword of a power of $b$, say $b^n$.
That is, $va^N$ is an initial segment of $b^n$ for some $v$.
For each $i=0,\dots , N$ there exists an element $e_i\in F$ with
$l(e_i)\le l(b)$ such that $va^ie_i$ is a power of $b$, that is
$va^ie_i=b^{j_i}$. All $N+1$ elements $e_0,\dots e_N$ have length at
most $l(b)$. Therefore by the choice of $N$ there are $i_1, i_2$,
$0\le i_1<i_2\le N$ such that $e_{i_1}=e_{i_2}=e$. Then
$va^{i_1}e=b^{k}$, $va^{i_2}e=b^s$ for some $k,s$. Therefore
$e^{-1}a^{i_2-i_1}e=b^{s-k}$ with $i_2-i_1>0$. This contradicts our
assumption that no power of $a$ is conjugate to a power of $b$.
\smallskip
Thus $a^N$ is never a subword of a power of $b$. By symmetry there is
$N_1>0$ such that $b^{N_1}$ is never a subword of a power of $a$.
Then $M=\max\{N,N_1\}$ satisfies the requirements of Lemma 6.8.
\enddemo

\proclaim{Lemma 6.9} Let $F=F(X)$ be a free non-Abelian group with basis
$X$.  Let $a$ be a nontrivial cyclically reduced word in $X$.
Then $a$ is not a subword of $a^{-n}$ for $n>0$ and $a^{-1}$ is not a
subword of $a^n$ for $n>0$.
\endproclaim
\demo{Proof}
Suppose that $a$ is a subword of $a^{-n}$ where $n>0$. Since
$l(a)=l(a^{-1})$ this implies that $a$ is in fact a cyclic permutation
of $a^{-1}$, that is $z^{-1} a z=a^{-1}$ for some $z\in F$. Clearly
$z$ does not commute with $a$ since $a$ has infinite order. Therefore
the subgroup $H=gp(a,z)$ of $F$ is free of rank 2 with basis $a, z$.
However in the free group on $a,z$ we obviously have $z^{-1} a z\ne
a^{-1}$ which gives us a contradiction.
Lemma 6.9 is proved.
\enddemo

\proclaim {Lemma 6.10} Let $F=F(X)$ be a free non-Abelian group of finite
rank with basis $X$. Let $a$ and $s$ be nontrivial cyclically reduced
words in $X$. Let $N$ be the number of distinct elements of
$F$ of length at most $l(a)$. 
Suppose that $s^N$ is a subword of $a^n$ for some $n$. Then a
nontrivial power of $s$ is conjugate to a power of $a$.
\endproclaim
\demo{Proof}
The proof is exactly the same as the proof of Lemma 6.8 given above.
\enddemo

\proclaim{Lemma 6.11} Let $F=F(X)$ be a free group of finite rank and let $a$ be a nontrivial root-free cyclically reduced word in $X$. Let $N$ be the number of elements in $F(X)$ of length at most $l(a)$. Suppose $s$ is a nontrivial cyclically reduced word such that some positive power $a^n$ of $a$ has initial segment $s^{N+1}$. Then $s=a^j$ for some $j>0$.
\endproclaim 
\demo{Proof} 
Suppose $s$ and $a$ are as above and so $s^{N+1}$ is an initial segment of $a^n$. Then for every $i=1,2,\dots N$ there is a word $e_i$ with $l(e_i)\le l(a)$ such that $s^ie_i=a^{n_i}$ for some $0\le n_i\le n$.
by the choice of $N$ this implies that there are some $i_1, i_2$, $1\le i_1<i_2\le N$ such that $e_{i_1}=e_{i_2}=e$. Then $e=s^{-i_1}a^{n_{i_1}}$ and $e a^{n_{i_2}-n_{i_1}} e^{-1}=s^{i_2-i_1}$.
Hence

$$s^{-i_1}a^{n_{i_1}} a^{n_{i_2}-n_{i_1}} a^{-n_{i_1}}s^{i_1} =s^{i_2-i_1}$$
and
$$a^{n_{i_2}-n_{i_1}}=s^{i_2-i_1}.$$
Since $a$ is a root-free element and the subgroup $gp(a)$ is malnormal in $F$, this implies that $s=a^j$ for some $j>0$.
Lemma 6.11 is proved.

\enddemo
\proclaim{Lemma 6.12} Let $F=F(X)$ be a free group of finite rank and let $a$ be a nontrivial root-free cyclically reduced word in $X$. Let $N$ be the number of elements in $F(X)$ of length at most $l(a)$.
Suppose $s$ is a nontrivial freely reduced word in $X$ such that the product $s\cdot a$ is freely reduced and such that $s a^n$ has initial segment $a^n$ for some $n$ with $nl(a)\ge (N+1) l(s)$. Then $s\equiv a^j$ for some $j>0$.
\endproclaim

\demo{Proof}
Lemma 6.7 implies that $s$ is cyclically reduced and that $a^n\equiv s^k\alpha$ with $l(\alpha)<l(s)$. Thus
$$(k+1)l(s)> l(a^n)=nl(a)\ge (N+1)l(s) \text{  by the choice of } n$$
and therefore $k\ge N+1$. Since $a^n$ has initial segment $s^k$, Lemma 6.11 implies that $s=a^j$ for some $j>0$.
Lemma 6.12 is proved.
\enddemo

\proclaim{Proposition 6.13} Let $F=F(X)$ be a non-abelian free group of finite rank and let $a\in F$ be a non-trivial cyclically reduced root-free word in $X$.
Then there exists a positive integer $K$ such that the following holds.

Let $n\ge K$ and let $H=<a^n>=gp(a^n)$. Let $g\in F$ be a shortest element in the double coset $HgH$ and such that the element $g\cdot a^n$, when freely reduced over $X$, has $a^{[5n/8]}$ as its initial segment. 

Then $g\in gp(a)=<a>$.
\endproclaim 

\demo{Proof}
Let $N$ be the number of elements of $F(X)$ of length at most $l(a)$. Put $K=100(N+1)$ and assume that $n\ge K$.

Let $g\equiv x y^{-1}$ and $a^n\equiv y z$ where $y$ is the maximal initial segment of $a^n$ that is freely cancelled in $g \cdot a^n$.
Note that $l(y)\le (1/2) l(a^n)$ because $g$ is shortest in $gH$.
\smallskip
We have that $g a^n\equiv x z$ has $a^{[5n/8]}$ as its initial segment.
Thus $x$ has length $l(x)\le (1/2)l(a^n)$ because $g\equiv x y^{-1}$ is shortest in $Hg$. This implies, in particular, that $x$ is an initial segment of $a^{[5n/8]}$.
We also conclude that $l(g)\le l(x)+l(y)\le l(a^n)$.
\smallskip

Since $y$ is an initial segment of $a^n$ with $l(y)\le l(a^n)/2$, it has the form
$$y\equiv a^k a', \quad\text{ where }\quad a^n\equiv a^ka'a'' a^{n-k-1},\quad a\equiv a'a'',\quad\text{ and }k\le n/2.$$

Similarly, since $x$ is an initial segment of $a^{[5n/8]}$ of length at most $l(a^n)/2$, we have

$$x\equiv a^s a_1, \quad\text{ where }\quad a^n\equiv a^sa_1a_2 a^{n-s-1},\quad a\equiv a_1a_2,\quad\text{ and }s\le n/2.$$

Thus $g\equiv x y^{-1}\equiv a^s a_1 (a')^{-1} a^{-k}$ and
$$g a^n=a^s a_1 (a')^{-1} a^{-k}\cdot a^ka'a'' a^{n-k-1}\equiv a^s a_1a'' a^{n-k-1}.$$

Since $g a^n$ starts with $a^{[5n/8]}$ and $s\le n/2$ this implies that
$a_1a'' a^{n-k-1}$ has initial segment $a^{[n/8]}$. Since $k\le n/2$ this implies that $a_1a'' a^{[n/8]}$ has initial segment $a^{[n/8]}$.
Note that $a_1$ and $a''$ are segments of $a$, so that $l(a_1a'')\le 2l(a)$.
By the choice of $n\ge K=100(N+1)$ we have
$$[n/8] l(a)> (N+1) (2l(a))\ge (N+1)l(a_1a'').$$

Therefore Lemma 6.12 implies that $a_1a''\equiv a^j$ for some $0\le j\le 2$.
\smallskip

Case 1. Suppose that $j=2$. Then $l(a_1)=l(a'')=l(a)$ and therefore $a_1=a''=a$, $a_2=a'=1$. Hence $g=a^s a_1 (a')^{-1} a^{-k}=a^s a a^{-k}=a^{s-k+1}$ and $g\in gp(a)$ as required.
\smallskip

Case 2. Suppose that $j=1$ and $a_1a''\equiv a$. Since $a_1a_2\equiv a'a''\equiv a$, this implies that $a'=a_1, a''=a_2$. Hence $g=a^s a_1 (a')^{-1} a^{-k}=a^s a_1 (a_1)^{-1} a^{-k}=a^{s-k}\in gp(a)$ as required.
\smallskip

Case 3. Suppose that $j=0$ and $a_1=a''=1$. Then $a_2=a'=a$ and hence $g=a^s a_1 (a')^{-1} a^{-k}= a^s a^{-1} a^{-k}=a^{s-k-1}\in gp(a)$ as required.

This completes the proof of Proposition 6.13.

\enddemo

\proclaim{Proposition 6.14} Let $F=F(X)$ be a non-abelian free group of finite rank. Let $a$ and $b$ be nontrivial cyclically reduced elements of $F$ such that no nontrivial power of $a$ is conjugate in $F$ to a power of $b$.Then there exists an integer $N>1$ such that for any $n\ge N, m\ge N$, $A=a^n$, $B=B^m$, $H=gp(a^n, b^m)=gp(A,B)$ the following holds.

If $g\in F$ is shortest in the double coset $HgH$ and if $g H g^{-1}\cap H\ne 1$ then $g\in gp(a)\cup gp(b)$.
\endproclaim 
\demo{Proof}
Let $K_0$ be such that
\roster
\item the length of the maximal initial segment of $b^j$ which may be freely cancelled in $a^i \cdot b^j$, $i,j\in{\Bbb Z}$ is less than $K_0$;
\item $K_0l(a)\ge l(b)$ and $K_0l(b)\ge l(a)$
\item $a^{K_0}$ is not a subword of a power of $b$ and $b^{K_0}$ is not a subword of a power of $a$.
\endroster
Let $K_a$ and $K_b$ be the constants provided by Proposition 1 for $a$ and $b$ respectively. Put $N=1000 K_0 \max\{K_a,K_b\} l(a) l(b)$.
Let $n,m\ge N$, $A=a^n$, $B=b^m$, $H=gp(A,B)$. Let $g$ be as in Proposition 2.
\smallskip

Notice that the choice of $n,m\ge N$ implies that the length of the maximal initial segment of $B^{\pm 1}$ that may be freely cancelled in the various products $A^{\pm 1} B^{\pm 1}$ is less than $(1/1000)\min\{l(A), l(B)\}$.

Since $gHg^{-1}\cap H\ne 1$ we have that
$g W(A,B) g^{-1}=W_1(A,B)$ for some freely reduced words $W,W_1$ in $A$, $B$.
By taking powers, if necessary, we may assume that the length of $W$ (as an $A,B$-word) is greater than $1000 l(g)$.
We will rewrite the equality above as

$$ g W(A,B)= W_1(A,B) g=f.$$

Without loss of generality we may assume that $W$ starts with $A$.
There are several cases to consider.
\smallskip

Case 1. The word $W_1(A,B)$ starts with $A$. Then $f$ has initial segment $a^{[999/1000 n]}$. Let $g\equiv x y^{-1}$ where $y$ is the maximal initial segment of $A$ freely cancelled in $g\cdot A$. Since $g$ is shortest in $gp(A,B)\cdot g \cdot gp(A,B)$, we have $l(g)\le (n/2)l(a)=1/2 l(A)$.

Subcase 1.a. Suppose that $l(x)> 251/1000 n l(a)$. Then $A=a^n=yzq$ where $q$ is the terminal segment of $A$ that is cancelled in the product of $A$ with the second letter of $W$. As we noted before, we have $l(q)< (1/1000)l(A)$ and $l(y)\le (1/2)l(A)$. Thus $l(z)> 499/1000 l(A)$, $l(x)>251/1000 l(A)$ and so $l(xz)> 750/1000 l(A)=3/4 l(A)$.

Obviously, $xz$ is an initial segment of $gW(A,B)=f$ and of $gA$ of length greater than $[3/4 n]l(a)$. We have also seen that $g W(A,B)=f$ has initial segment $a^{[999/1000 n]}$.
Thus $gA=g a^n$ has initial segment $a^{[3/4n]}$. Since $g$ is shortest in $gp(A,B)\cdot g \cdot gp(A,B)$, the element $g$ is also shortest in $gp(A)\cdot g \cdot gp(A)$. Therefore by Proposition 6.13 we have $g\in gp(a)\subset gp(a)\cup gp(b)$ as required.
\smallskip

Subcase 1.b. Suppose that $l(x)\le 251/1000 l(A)$.

If $l(y)> 251/1000 l(A)$ then we may replace the identity
$g W(A,B)= W_1(A,B) g$ by $g^{-1} W_1(A,B)= W(A,B) g^{-1}$. Since $g^{-1}$ is obviously also shortest in $Hg^{-1}H$, we observe that this case is symmetric to Subcase 1.a and so $g^{-1}\in gp(a), g\in gp(a)$.
\smallskip
Suppose, on the contrary, that $l(y)\le 251/1000 l(A)$.
Recall that $A=yzq$ where $l(q)<1/1000 l(A)$. Therefore $l(z)> 748/1000 l(A)$ and so $l(xz)> 748/100 l(A)$.
We see that both $f=gW=W_1g$ and $gA$ have initial segment $xz$ of length greater than $748/1000 l(A)$. On the other hand $gW=f$ has initial segment $a^{[999/1000 n]}$. Therefore $gA=ga^n$ has initial segment $a^{[748/1000n]}$ and so $a^{[5/8n]}$.  Since $g$ is shortest in $gp(A,B)\cdot g \cdot gp(A,B)$, the element $g$ is also shortest in $gp(A)\cdot  g\cdot gp(A)$. Therefore by Proposition 6.13 we have $g\in gp(a)\subset gp(a)\cup gp(b)$ as required.
\smallskip

Case 2. Suppose that the first letter of $W_1(A,B)$ is $B$.

Recall that $g=x y^{-1}$, $A=yzq$ where
$$l(q)<(1/1000) l(A), l(y)\le 1/2 l(A), l(z)> (499/1000) l(A).$$
Also $gW=f$ has initial segment $xz$.
On the other hand $f=W_1g$ has initial segment $b^{[999/1000n]}$ because $W_1(A,B)$ starts with $B$. Hence $l(x)\le 1/2l(B)$ because otherwise $g\equiv x y^{-1}$ would not have been shortest in $gp(B)\cdot g \cdot gp(B)$.
Thus $x$ is an initial segment of $B$ with $l(x)\le 1/2l(B)$.
Since $l(z)>499/1000l(A)$, the choice of $n\ge N$ implies that $l(z)>l(b)(N_0+1)$. Let $z'$ be the initial segment of $z$ of length $l(b)(N_0+1)$.

Then $xz'$ and $b^{[999/1000m]}$ are initial segment of $f$ where $l(x)\le 1/2ml(b)$. Therefore $z'$ is a subword of $b^{[999/1000m]}$ and of $a^n$ of length $l(b)(N_0+1)$. This implies that $b^{N_0}$ is a subword of $a^n$ which is impossible by the choice of $N_0$.
\smallskip

Case 3. Suppose the first letter of $W_1(A,B)$ is $B^{-1}$. This case is completely analogous to Case 2.
\smallskip

Case 4. Suppose the first letter of $W_1(A,B)$ is $A^{-1}$. Then, once again,
$g\equiv x y^{-1}$, $A=yzq$ where
$$l(q)<1/1000 l(A), l(y)\le 1/2 l(A), l(z)> 499/1000 l(A).$$
Also $gW=f$ has initial segment $xz$.
On the other hand $f=W_1g$ starts with $a^{-[999/1000n]}$. Hence $l(x)\le 1/2l(A^{-1})=l(A)$ since $g=xy^{-1}$ is shortest in $gp(A)\cdot g \cdot gp(A)$.
Let $\hat z$ be an initial segment of $z$ of length $2l(a)$ (and so $\hat z$ contains the subword $a$). Then $\hat z$ is a subword of $a^{-[999/1000 n]}$ and so $a$ is a subword of $a^{-n}$, which is impossible by Lemma 6.9.

This completes the proof of Proposition 6.14.

\enddemo

\proclaim{Proposition 6.15} Let $F=F(X)$ be a non-abelian free group of finite rank. Let $a,b\in F$ be non-trivial root-free elements such that no non-trivial power of $a$ is conjugate to a power of $b$. Then there exists an integer $N>1$ with the following properties. 
\smallskip
Suppose $m,n\ge N$ and let $D=gp(a^nb^{3m}a^n, b^m a^{3n} b^m)$. Then $D$ is malnormal in $F$.
\endproclaim
\demo{Proof}
Let $N$ be the constant provided by Proposition 6.14. Suppose that $D$ is not malnormal in $F$ and there is $g\in F-D$ such that $g D g^{-1}\cap D\ne 1$.
Put $H=gp(a^n,b^m)$, so that $D\le H$. Let $g'$ be the shortest element in $HgH$ so that $g=h_1 g' h_2$ for some $h_1,h_2\in H$. The fact that $g D g^{-1}\cap D\ne 1$ implies $g H g^{-1}\cap H\ne 1$. Thus
$$\align
h_1 g' h_2 &H h_2^{-1} (g')^{-1}h_1^{-1}\cap H\ne 1\\
h_1g' H (g')^{-1}h_1^{-1}\cap H\ne 1\quad\text{ and so }\\
g' H (g')^{-1}\cap H\ne 1\\
\endalign$$

By Proposition 6.14 this means that $g'\in gp(a)\cup gp(b)$ and so $g'\in gp(a,b)=L$. Since $H=gp(a^n,b^m)\le gp(a,b)=L$ and $g=h_1 g' h_2$ with $h_1,h_2\in H$, we have $g\in L$. 
Thus $D\le L$, $gDg^{-1}\cap D\ne 1$ and $g\in L$.
However, $D$ is malnormal in $L$ by Lemma 6.6. Therefore $g\in D$ which contradicts our assumption that $g\in F-D$.
Proposition 6.15 is proved.
\enddemo

\proclaim{Theorem 6.16 (Theorem C)} Let $F=F(X)$ be a non-Abelian free group of finite
rank with basis $X$.  Let
$H_1, \dots H_k$ be finitely generated subgroups of $F$ of infinite
index. Then there exist elements $a,b\in F$ such that
\roster
\item the subgroup $H=gp(a,b)$ is a free group of rank two which is
malnormal in $F$;
\item no nontrivial element of $H$ is conjugate in $F$ to an element
of $H_1\cup\dots \cup H_k$.
\endroster
\endproclaim

\demo{Proof}
By Lemma 6.1 there exists $y\in F$ such that no nontrivial power of
$y$ is conjugate in $F$ to an element of $H_1\cup\dots \cup H_k$.
Take $t$ to be any cyclically reduced word in $X$ such that no power of $t$
is conjugate in $F$ to a power of $y$. (Since $|X|>2$ we may even
choose $y\in X$.) 
Put $D_n=gp(y^nt^{3n}y^n, t^n y^{3n} t^n)\le F$ where $n>0$.
By Lemma 6.3 there exists $M_1>0$ such that if $n\ge M_1$ then no
nontrivial element of $D_n$ is conjugate in $F$ to an element of
$H_1\cup\dots \cup H_k$.
By Proposition 6.15 there is $M_2>0$ such that in $n\ge M_2$ then
$D_n$ is malnormal in $F$. Put $M=\max \{M_1, M_2\}$.
Then for any $n\ge M$ the subgroup $D_n$ of $F$ satisfies the
requirements of Theorem 6.16.
\enddemo

\head 7.Quasiconvex subgroups and malnormality \endhead

In this section we will show that if $G$ is a torsion-free hyperbolic
group and $H$ is a non-elementary subgroup of $G$ then $H$ contains a
free subgroup $F$ of rank two such that $F$ is malnormal and
quasiconvex in $G$.
\smallskip
First, we need to accumulate a certain amount of information about the
properties of quasiconvex subgroups.
\proclaim{Lemma 7.1} Let $G$ be a word hyperbolic group with a finite
generating set $X$ and word metric $d_X$. Let $H$ be a
quasiconvex subgroup of $G$. Then there exists $\lambda >0$ with the
following properties.

Suppose that $g\in G$ is such that $g\not\in H$ and that $g$ is
shortest with respect to $d_X$ in the double coset class $HgH$.
Choose $w$ to be an $X$-geodesic representative of $g$. Then for any
$X$-geodesic words $v_1$ and $v_2$ representing elements of $H$ the
word $v_1wv_2$ is $(\lambda, \lambda)$-quasigeodesic with respect to
$d_X$.
\endproclaim
\demo{Proof}

The statement of Lemma 7.1 easily follows from the proof of Lemma 4.5 in
\cite{BGSS} and we omit the details.
\enddemo
\remark{Remark} Suppose that $G$ is a group and $M$ is a subgroup of
$G$. Let $z\in G$ and $z'=m_1zm_2$ for some $m_1,m_2\in M$.
Then $zMz^{-1}\cap M\ne 1$ if and only if ${z'}M{z'}^{-1}\cap M\ne 1$.
\endremark

The following lemma shows that if $H$ is a quasiconvex subgroup of a
torsion-free hyperbolic group $G$ then there are at most finitely many
double cosets $HzH$ that ``violate'' malnormality of $H$.

\proclaim{Lemma 7.2} Let $G$ be a torsion-free word hyperbolic group
with a finite generating set $X$ and word metric $d_X$. Let $H$ be a
nontrivial quasiconvex subgroup of $G$. Then there are only finitely
many double cosets $HzH$ of elements $z\in G$, $z\not\in H$ such that
$z^{-1}Hz\cap H\ne 1$.
\endproclaim
\demo{Proof}
Let ${\Cal H}=\{h_1, \dots h_k\}$ be a finite generating set of
$H$. For every $i=1,\dots k$ pick a $d_X$-geodesic word $v_i$
representing $h_i$. Since $H$ is quasiconvex in $G$, there is
${\lambda}_1 >0$ such that if $W(h_1,\dots ,h_k)$ is an ${\Cal H}$-geodesic word then $W(v_1,\dots , v_k)$ is a $(\lambda_1,
\lambda_1)$-quasigeodesic with respect to $d_X$.
\smallskip
Recall also that by Lemma 7.1 there is $\lambda >0$ such that if $v$ is
an $X$-geodesic representative of an element $z\in G$, $z\not\in L$
which is $d_X$-shortest in $HzH$,
then for any $X$-geodesic words $q_1, q_2$ representing elements of
$H$ the word $q_1vq_2$ is $(\lambda, \lambda)$-quasigeodesic with
respect to $d_X$.
\smallskip

Put $\lambda'=\max\{\lambda, \lambda_1\}$. Let $\epsilon >1$ be such
that any two $(\lambda', \lambda')$-quasigeodesic paths in the Cayley
graph $\Gamma(G,X)$ with common endpoints are $\epsilon$-Hausdorff
close.
Let $K$ be the maximal length of the words $v_1,\dots v_k$.
\smallskip

Suppose now that there are infinitely many double cosets $HzH$ such
that $$z H z^{-1}\cap H\ne 1.$$ Then there exists $z\in G$, $z\not \in H$
such that $z H z^{-1}\cap H\ne 1$, $z$ is shortest in $HzH$ and
$l_X(z)>2\epsilon +K+2$. Let $Z$ be a $d_X$-geodesic representative of
$z$.
\smallskip
Then $zhz^{-1}=h'$  and $zh=h'z$  for some $h,h'\in H$. Replacing $h$ and $h'$
by some high powers of themselves we can assume that $l_X(h')>3\epsilon+K+3$.
\smallskip
Choose a $d_X$-geodesic
representative $V$ for $h$ and a $d_X$-geodesic representative $V'$
for $h'$. 
Then both words $ZV$ and $V'Z$ are $(\lambda, \lambda)$-quasigeodesic
with respect to $d_X$ representing the element $zh=h'z$.  Recall that $l_X(z)>2\epsilon +K+2$. Take $Z_1$
to be the initial segment of $Z$ of length $2\epsilon +K+2$ so that
$Z\equiv Z_1Z_2$ (see Figure 7.1).
\midinsert
\epsfxsize=12cm\epsfbox{f7.1.eps}
\captionwidth{6cm}
\botcaption{Figure 7.1}
\endcaption
\endinsert
By the choice of $\epsilon$ there is an initial
segment $p$ of the word $V'Z$ such that $d_X(\overline{Z_1},\overline
p)\le \epsilon+1$. Since $l(Z_1)=2\epsilon +K+2$, this implies that
$l_X(\overline p)\le 2\epsilon +K+2+\epsilon+1=3\epsilon+K+3$. Note
that $h'$ was chosen so that $l_X(h')>3\epsilon+K+3$. Therefore the
initial segment $p$ of $V'Z$ is in fact an initial segment of $V'$.
\medskip
Note that $zh=h'z$.
Choose  
a $d_{\Cal H}$ geodesic representative
$W(h_1,\dots h_k)$ of $h'$. Denote the $X$-word $W(v_1,\dots , v_k)$
by $U$. Then both $U$ and $V'$ represent the element
$h'$ of $G$. Moreover, both words $U$ and  $V'$ are
$(\lambda', \lambda')$-quasigeodesic with respect to $d_X$ and so are
$\epsilon$-Hausdorff close in $\Gamma(G,X)$. Therefore 
by the choice of $\epsilon$ and $K$ there is an initial segment
$W_1(h_1,\dots , h_k)$ of $W(h_1,\dots , h_k)$ such that
$d_X(\overline{p}, \overline{W_1})\le \epsilon+K$. 
Therefore $d_X(\overline{Z_1}, \overline{W_1})\le
d_X(\overline{Z_1},\overline p)+d_X(\overline{p}, \overline{W_1})\le 
\epsilon+1+\epsilon+K=2\epsilon +K+1$. 
\smallskip
Denote $q=(\overline{W_1})^{-1}\overline{Z_1}\in G$. Then $$l_X(q\overline{Z_2})\le 2\epsilon+K+1+l(Z_2)<2\epsilon+K+2+l(Z_2)=l(Z_1)+l(Z_2)=l(Z)=l_X(z).$$

On the other hand $\overline{W_1}\cdot q\overline{Z_2}=z$ and
$\overline {W_1}\in H$. Thus $Hz=Hq\overline{Z_2}$ and so
$HzH=Hq\overline{Z_2}H$. But we have established that
$l_X(q\overline{Z_2})<l_X(z)$. This gives us a contradiction with the
choice of $z$ as a shortest element in $HzH$. Lemma 7.2 is proved.
\enddemo   

\proclaim{Lemma 7.3} Let $G$ be a torsion-free word hyperbolic group and
let $H$ be a nontrivial quasiconvex subgroup of $G$. Let $g\in G$,
$g\not\in H$  be such that $g^iHg^{-i}\cap H\ne 1$ for every $i>0$.
Then for some $n>0$ we have $g^n\in H$.
\endproclaim
\demo{Proof}
This statement is an immediate corollary of the main result of the
paper \cite{GMSR} by R.Gitik, M.Mitra, M.Sageev and E.Rips.
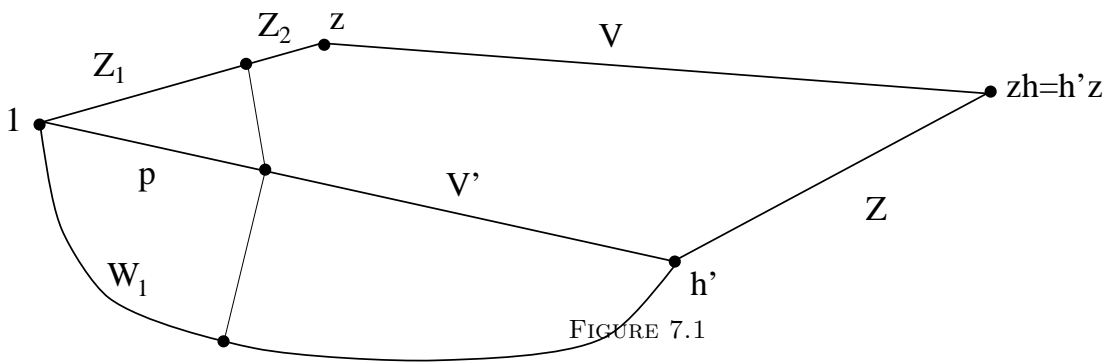
\enddemo
   
The following statement is essentially due to M.Gromov \cite{Gr, Theorem 5.3.E}. For a
careful argument the reader is referred to Lemma 1.1 and Lemma 1.2 in the paper of T.Delzant \cite{D}.
\proclaim{Lemma 7.4} Let $G$ be a word hyperbolic group. Let $a, b\in G$
be such that no nontrivial power of $a$ is equal to a power of
$b$. Then there is $n>0$ such that the subgroup $H=gp(a^n, b^n)$ is
free of rank two and is quasiconvex in $G$.
\endproclaim  

\proclaim{Lemma 7.5} Let $G$ be a torsion-free word hyperbolic
group. Let $H=F(a,b)$ be a quasiconvex subgroup of $G$ which is free
group with basis $a,b$. Put $V_G(H)$ to be the virtual normalizer of $H$ in $G$,
that is
$$V_G(H)=\{g\in G \mid \ |H: H\cap gHg^{-1}|<\infty, |gHg^{-1}:H\cap gHg^{-1}|<\infty\}.$$

Then $H=V_G(H)$.
\endproclaim

\demo{Proof}
Since $H$ is infinite and quasiconvex in $G$, the subgroup $H$ has
finite index in its virtual normalizer $V_G(H)$ (see Proposition 3.4). Also, $G$ is
torsion-free and $H$ is a free group. Thus $V_G(H)$ is a torsion-free
group that has a free subgroup of finite index. By the theorem of
J.Stallings \cite{Sta} this implies that $V_G(H)$ is itself a free group. 
Suppose that $V_G(H)\ne H$, that is the index of $H$ in $V_G(H)$ is
greater than 1. By the Theorem ? of \cite{LS} this implies that the rank of
$V_G(H)$ is strictly less than the rank of $H$. However the rank of
$H$ is equal to two and $H\le V_G(H)$ which gives us a contradiction.
Thus $H=V_G(H)$ and Lemma 7.5 is proved.
\enddemo

\proclaim{Lemma 7.6} Let $G$ be a group and let $H$ be a subgroup of
$G$. Let $g\in G$ and suppose that the subgroup $K_1=gHg^{-1}\cap H$ has
finite index in $H$ and infinite index in $gHg^{-1}$. Then for any
integer $n>0$ the subgroup $K_n=g^nHg^{-n}\cap H$ has finite index in $H$
and infinite index in $g^nHg^{-n}$.
\endproclaim
\demo{Proof}
We will prove the statement by induction on $n$. When $n=1$, the
statement is obvious. Suppose now $n>1$ and that Lemma 7.6 has been
established for all smaller values.

Thus by the inductive hypothesis the subgroup
$K_{n-1}=g^{n-1}Hg^{-(n-1)}\cap H$ has finite index in $H$ and
infinite index in $g^{n-1}Hg^{-(n-1)}$. Conjugating everything by $g$
we conclude that $gK_{n-1}g^{-1}=g^nHg^{-n}\cap gHg^{-1}$ has finite
index in $gHg^{-1}$ and infinite index in $g^nHg^{-n}$.
\smallskip
We know that the subgroup $K_1=gHg^{-1}\cap H$ has finite index in
$H$. Also the intersection of $gHg^{-1}$ and $g^nHg^{-n}$ has infinite
index in $g^nHg^{-n}$. Therefore the intersection of any subgroup of
$gHg^{-1}$ and $g^nHg^{-n}$ has infinite
index in $g^nHg^{-n}$. In particular $K_1\cap g^nHg^{-n}$ has infinite
index in $g^nHg^{-n}$. Since $K_1$ is a subgroup of finite index in
$H$, this implies that $K_n=H\cap g^nHg^{-n}$ has infinite
index in $g^nHg^{-n}$.
\smallskip
On the other hand we already observed that
$L_1=gK_{n-1}g^{-1}=g^nHg^{-n}\cap gHg^{-1}$ has finite index in
$gHg^{-1}$. Recall that $K_1$ is a subgroup of $gHg^{-1}$ and therefore $K_1\cap
L_1$ has finite index in $K_1$. Since $K_1=gHg^{-1}\cap H$ is of
finite index in $H$, this implies that $K_1\cap L_1$ 
 is of finite index in $H$. However $K_1\cap L_1\le L_1\le g^nHg^{-n}$
and therefore $g^nHg^{-n}\cap H$ is of finite index in $H$. This
completes the inductive step and the proof of Lemma 7.6.
\enddemo

\proclaim{Theorem 7.7} Let $G$ be a torsion-free word hyperbolic group. Let $\Gamma$
be a non-elementary subgroup of $G$. Then there is a subgroup $F\le
\Gamma$ such that $F$ is a free group of rank 2 , $F$ is malnormal in
$G$ and $F$ is quasiconvex in $G$.
\endproclaim

\demo{Proof}

All nontrivial abelian subgroups of torsion-free hyperbolic groups are
infinite cyclic \cite{KM} and therefore elementary. Thus $\Gamma$ is
non-abelian. Therefore $\Gamma$ contains two non-commuting elements
$a_0$, $b_0$. Notice that no nontrivial power of $a_0$ is equal to a
power of $b_0$. Indeed, if $a_0^n=b_0^m$ for some $n>0$ then the
centralizer $C_G(a_0^n)$ of $a_0^n$ contains both $a_0$ and
$b_0$. However, centralizers of nontrivial elements in torsion-free
hyperbolic groups are infinite cyclic [] and therefore $C_G(a_0^n)$ is
cyclic and abelian. This implies that $a_0$ and $b_0$ commute which contradicts
our assumptions.

Therefore by Lemma 7.4 there exists $n>0$ such that the subgroup
$K=gp(a_0^n, b_0^n)$ is free of rank two and quasiconvex in
$G$. Denote $a=a_0^n$, $b=b_0^n$ so that $K=gp(a,b)=F(a,b)$ is a free
group on $a,b$.

By Lemma 7.5 the subgroup $K$ coincides with its virtual normalizer. If
$K$ is malnormal in $G$ then $K$ satisfies all the requirements of
Theorem 7.7 and there is nothing to prove. Suppose now that $K$ is not
malnormal in $G$. By Lemma 7.2 there exist only finitely many elements
double cosets $KzK$ where $z\in G, z\not \in K$ is such that
$zKz^{-1}\cap K\ne 1$. Let us denote these cosets $Kz_1K, \dots
Kz_sK$.

Observe that if $z\not \in K$ and $K\cap zKz^{-1}\ne 1$ then the
subgroup $K_1=K\cap zKz^{-1}$ has infinite index in both $K$ and
$zKz^{-1}$. Indeed, if both these indices are finite then $z\in
V_G(K)=K$.
Suppose that one of these indices is finite and the other is
infinite. Without loss of generality (replace $z$ by $z^{-1}$ if necessary) we may assume that $K_1$ has
finite index in $K$ and infinite index in $zKz^{-1}$. Then by Lemma 7.6
 for every $i>0$ the subgroup $K_i=K\cap z^i K z^{-i}$ has
finite index in $K$ and infinite index in $z^i K z^{-i}$. In
particular $K_i\ne 1$.
By Lemma 7.3 this implies that $z^p\in K$ for some $p>0$.  However this
means that $K_p=K=z^p K z^{-p}=K\cap z^p K z^{-p}$ has finite index
(in fact index one) in $z^p K z^{-p}$. This gives us a contradiction.
Thus we have established that for every $z\in G, z\not \in K$ such that
$K\cap zKz^{-1}\ne 1$, the subgroup $K\cap zKz^{-1}$ has infinite
index in both $K$ and  $zKz^{-1}$.
Notice also that since $K$ and $zKz^{-1}$ are quasiconvex in $G$,
their intersection $K\cap zKz^{-1}$ is also quasiconvex and so
finitely generated. 
\medskip
Thus all the subgroups $H_i=z_iKz_i^{-1}\cap K$ are finitely generated
and of infinite index in $K$. By Theorem 6.16 there exists a subgroup
$F$ of $K$ such that $F$ is free of rank two, $F$ is malnormal in $K$
and no nontrivial element of $F$ is conjugate in $K$ to an element of
$H_1\cup \dots \cup H_s$. 
\medskip
We claim that $F$ is malnormal in $G$. Indeed, suppose not. Then there
exist $f_1,f_2\in F$, $f_1\ne 1$ and $g\in G$, $g\not\in F$ such that
$g f_1 g^{-1}=f_2$. Since $F$ is malnormal in $K$, this implies
$g\not\in K$. Therefore $g$ has the form $g=k_1 z_i k_2$ for some
$k_1,k_2\in K$, $1\le i\le s$. Hence

$$ k_1 z_i k_2 f_1 k_2^{-1} z_i^{-1} k_1^{-1}=f_2$$
and

$$ z_i (k_2 f_1 k_2^{-1}) z_i^{-1} = k_1^{-1} f_2 k_1.$$

Therefore $k_1^{-1} f_2 k_1\in z_iKz_i^{-1}\cap K=H_i$. However
$f_2\ne 1, f_2\in F$ and the subgroup  $F$ of $K$ was chosen so that
no nontrivial element of $F$ is conjugate in $K$ to an element of
$H_i$. This gives us a contradiction.
So $F$ is indeed malnormal in $G$.
\medskip

Any finitely generated subgroup in a finite rank free group is
quasiconvex \cite{Sho} and therefore $F$ is quasiconvex in $K$. Since $K$ is
quasiconvex in $G$, this implies that $F$ is quasiconvex in $G$ (see
Proposition 3.2(6)).
Thus $F\le K\le \Gamma$ is a free group of rank two that is malnormal
and quasiconvex in $G$. Theorem 7.7 is proved.
\enddemo

\enddocument